\documentclass[12pt,reqno]{amsart}
\usepackage{amsmath,amssymb}
\usepackage{graphicx,color}
\usepackage{amscd,amsthm,amsfonts,boxedminipage}
\usepackage{boxedminipage}
\usepackage{array}

\input xy
\xyoption{all}

\newcommand{\R}{{\mathbb{R}}}
\newcommand{\Z}{{\mathbb{Z}}}
\newcommand{\C}{{\mathbb{C}}}

\newcommand{\bea}{\begin{eqnarray}}
\newcommand{\eea}{\end{eqnarray}}

\newcommand{\bp}{\begin{pmatrix}}
\newcommand{\ep}{\end{pmatrix}}

\newcommand{\bps}{\begin{smallmatrix}}
\newcommand{\eps}{\end{smallmatrix}}

\newcommand{\ti}{\tilde}

\def \cE{{\mathcal E}}
\def \cF{{\mathcal F}}

\def \cI{{\mathcal I}}

\def \cO{{\mathcal O}}

\def \cU{\mathcal{U}}

\def \cV{{\mathcal V}}

\def \f{{\frak f}}

\def \sbar{\underline{s}}

\def \0{{\bf 0}}
\def \1{{\bf 1}}
\def \ii{{\bf i}}

\def \DG{\mathit{DG}}

\def \fpartial#1{\frac{\partial}{\partial {#1}}}
\def \fpart#1#2{\frac{\partial #1}{\partial #2}}

\def \bpartial{{\bar{\partial}}}

\def \ov#1{\frac{1}{#1}}
\def \({\left(}
\def \){\right)}

\def \ti#1{\tilde{#1}}

\def \zb{{\bar z}}
\def \f{{\frak f}}

\def \cO{{\mathcal O}}

\def \grad{\mathrm{grad}}
\def \graph{\mathrm{graph}}

\def \Int{\mathrm{Int}}

\def \Tri{\mathit{Tr}}

\def \Mo{\mathit{Mo}}

 \newtheorem{thm}{Theorem}[section]
 \newtheorem{lem}[thm]{Lemma}
 \newtheorem{prop}[thm]{Proposition}
 \newtheorem{cor}[thm]{Corollary}

\theoremstyle{definition}

 \newtheorem{rem}[thm]{Remark}

\begin{document}

\pagestyle{myheadings}

 \title{Homological mirror symmetry of {$\C P^n$} and their products via Morse homotopy}

\date{\today}

\noindent{
 }\\

\author{Masahiro Futaki}
\author{Hiroshige Kajiura}
\address{Graduate School of Science, Chiba University, 
  263-8522 Japan}
\email{futaki@math.s.chiba-u.ac.jp}
\email{kajiura@math.s.chiba-u.ac.jp}

\begin{abstract}
We propose a way of understanding homological mirror symmetry 
when a complex manifold is a smooth compact toric manifold.
So far, in many example,
the derived category $D^b(coh(X))$ of coherent sheaves
on a toric manifold $X$ is compared with 
the Fukaya-Seidel category of the Milnor fiber of the
corresponding Landau-Ginzburg potential. 
We instead consider the dual torus fibration $\pi:M\to B$ of
the complement of the toric divisors in $X$, 
where $\bar{B}$ is the dual polytope of the toric manifold $X$.
A natural formulation of homological mirror symmetry in this set-up
is to define $Fuk(\bar{M})$ a variant of the Fukaya category and show
the equivalence $D^b(coh(X))\simeq D^b(Fuk(\bar{M}))$. 
As an intermediate step, we construct 
the category $\Mo(P)$ of weighted Morse homotopy 
on $P:=\bar{B}$ as a natural generalization
of the weighted Fukaya-Oh category proposed by 
[M.~Kontsevich and Y.~Soibelman, 
In {\em Symplectic geometry and mirror symmetry}, page 203 (2001)]. 
We then show a full subcategory $\Mo_\cE(P)$ of $\Mo(P)$ 
generates $D^b(coh(X))$ 
for the cases $X$ is a complex projective space and their products. 
\end{abstract}


\renewcommand{\thefootnote}{\arabic{footnote}}
\setcounter{footnote}{0}

\maketitle


\markboth{MASAHIRO FUTAKI AND HIROSHIGE KAJIURA}{HMS OF {$\C P^N$} AND THEIR PRODUCTS VIA MORSE HOMOTOPY}

 \section{Introduction}

In this paper, we propose a way of understanding homological mirror symmetry 
for the case of smooth compact toric manifolds.
So far, in many studies,
the derived category $D^b(coh(X))$ of coherent sheaves
on a toric manifold $X$ is compared with 
the Fukaya-Seidel category of the Milnor fiber of the
corresponding Landau-Ginzburg potential. 
In this paper, we consider the dual torus fibration $\pi:M\to B$, 
in the sense of Strominger-Yau-Zaslow construction \cite{SYZ}, 
of the complement of the toric divisors in $X$,
where $P=\bar{B}$ is the dual polytope of the toric manifold $X$. 
Fukaya discussed the Calabi-Yau cases in \cite{fukaya:asymptotic} where 
the K\"ahler metric degenerates at singular fibers of the torus fibration. 
In our set-up, the K\"ahler metrics go to infinity at the boundaries $\partial(P)$. 
In \cite{fang08}, Fang discusses homological mirror symmetry 
of $\C P^n$ along this line.
There, he starts with considering line bundles on $\C P^n$ and the corresponding
Lagrangians in the mirror dual side. 
His idea of discussing the homological mirror symmetry 
is to consider the category of constructible sheaves as an intermediate step. 
We instead apply 
Kontsevich-Soibelman's approach \cite{KoSo:torus} to our case
and consider a category $\Mo(P)$ of Morse homotopy.
Namely, 
a natural formulation of homological mirror symmetry in our situation 
is to define a variant of the Fukaya category $Fuk(\bar{M})$ and show
the equivalence $D^b(coh(X))\simeq D^b(Fuk(\bar{M}))$. 
As an intermediate step, we construct 
the category $\Mo(P)$ of weighted Morse homotopy 
on $P$ as a natural generalization of the weighted Fukaya-Oh category
proposed in \cite{KoSo:torus}. 
We then show a full subcategory $\Mo_\cE(P)$ of $\Mo(P)$
generates $D^b(coh(X))$
for the case $X$ is a complex projective space and their products. 
For more general $X$, we may consider the Lagrangian sections 
discussed in \cite{chan09}, where
Chan discusses the correspondence between
holomorphic line bundles over projective toric manifolds and
Lagrangian sections in the mirror dual.
The relation of such an approach with Abouzaid's one \cite{abouzaid09}
is also mentioned there. 

While Abouzaid \cite{abouzaid09} employed the geometric perturbation technique 
to establish transversality 
and used the notion of $A_\infty$ pre-categories to avoid 
the self-intersection problem, 
we pick up finitely many Lagrangians explicitly and allowed clean intersections.
This simplifies the comparison of the symplectic and complex sides 
and makes the mirror functor more word-by-word.

This paper is organized as follows.
In section \ref{sec:SYZ}, we recall 
the SYZ torus fibration set-up \cite{SYZ} 
following \cite{LYZ,leung05}. 
There, a pair of dual torus fibrations $M\to B$ and $\check{M}\to B$
is defined. 
We use this set-up
by identifying $\check{M}$ with the complement of the toric divisors
in a toric manifold $X$. 
In section \ref{sec:Lag-Vect}, we recall the correspondence of 
Lagrangian sections of $M\to B$ and holomorphic line bundles on $\check{M}$, 
again,
following \cite{LYZ,leung05}. 
In the last subsection, we demonstrate a Lagrangian section to be
derived from the line bundle $\cO(k)$ restricted on the complement
of the toric divisors for $X=\C P^n$ 
by the correspondence above. 
In section \ref{sec:Cat}, we first recall DG-categories $\cF(M)$ and $\cV(\check{M})$ 
associated to
$M$ and $\check{M}$, respectively, in \cite{hk:fukayadeform}.
Kontsevich-Soibelman's approach for the homological mirror symmetry \cite{KoSo:torus}
proposes an intermediate category $\Mo(B)$ and the existence of
an $A_\infty$-equivalence
\[
 Fuk(M) \simeq \Mo(B) \overset{\sim}{\to} \cF(M) .
\]
This $\Mo(B)$ is called the weighted Fukaya-Oh category or
the category of weighted Morse homotopy on $B$. 
In subsection \ref{ssec:Mo}, we propose a modification $\Mo(P)$ of $\Mo(B)$ 
where $P=\bar{B}$ is the dual polytope of a smooth compact toric manifold
$X$. 
In the last section, we discuss the correspondence between
$D^b(coh(X))$ and $\Mo(P)$
when $X$ is a complex projective space or their products. 
In particular, we see that
we can take strongly exceptional collections $\cE$ of $D^b(coh(X))$
consisting of line bundles and the corresponding 
full subcategory $\Mo_\cE(P)$ of $\Mo(P)$
so that
\[
 \Tri(\Mo_\cE(P))\simeq D^b(coh(X)) 
\]
where $\Tri$ is the Bondal-Kapranov-Kontsevich construction \cite{BK:enhanced, kon94} 
of triangulated categories from $A_\infty$-categories. 

 \section{Toric manifolds and $T^n$-invariant manifolds}
\label{sec:SYZ}

 \subsection{Dual torus fibrations}
\label{ssec:dual}
 
In this subsection, we briefly review the SYZ torus 
fibration set-up \cite{SYZ}. 
For more details see \cite{LYZ,leung05}. 
We follow the convention of \cite{hk:fukayadeform}. 

Throughout this section, 
we consider an $n$-dimensional tropical Hessian manifold $B$ 
which we will define shortly below. 
Our goal of this subsection is then to construct torus fibrations 
$M$ and $\check{M}$, which are dual to each other, over the common base space $B$. 
A smooth manifold $B$ is called {\em affine} 
if $B$ has an open covering $\{U_\lambda\}_{\lambda\in\Lambda}$ 
such that the coordinate transformation is affine. 
This means that, for any $U_\lambda$ and $U_\mu$ such that 
$U_\lambda\cap U_\mu\ne\emptyset$, 
the coordinate systems 
$x_{(\lambda)}:=(x^1_{(\lambda)},\dots,x^n_{(\lambda)})^t$ 
and $x_{(\mu)}:=(x^1_{(\mu)},\dots,x^n_{(\mu)})^t$ are 
related to each other by 
\begin{equation}\label{coordinate-transform}
 x_{(\mu)}= \varphi_{\lambda\mu} x_{(\lambda)} + \psi_{\lambda\mu}
\end{equation}
with some $\varphi_{\lambda\mu}\in GL(n;\R)$ and $\psi_{\lambda\mu}\in\R^n$. 
If in particular $\varphi_{\lambda\mu}\in GL(n,\Z)$ for any 
$U_\lambda\cap U_\mu$, then $B$ is called {\em tropical affine}. 
(If in addition $\psi_{\lambda\mu}\in\Z^n$, $B$ is called {\em integral affine}. ) 

For simplicity, we take such an open covering $\{U_\lambda\}_{\lambda\in\Lambda}$ 
so that the open sets $U_\lambda$ and their intersections are all contractible. 
It is known 
that $B$ is an affine manifold iff 
the tangent bundle $TB$ is equipped with a torsion free flat connection. 
When $B$ is affine, then its tangent bundle $TB$ 
forms a complex manifold. 
This fact is clear as follows. 
For each open set $U=U_\lambda$, 
let us denote by $(x^1,\dots,x^n;y^1,\dots,y^n)$ 
the coordinates of $U\times\R^n \simeq TB|_{U}$ 
so that a point $\sum_{i=1}^n y^i\fpartial{x^i}|_x\in T_xB\subset TB$ 
corresponds to $(x^1,\dots,x^n;y^1,\dots,y^n)\in U\times\R^n$. 
We locally define the complex coordinate system by 
$z:=(z^1,\dots,z^n)^t$, where 
$z^i:=x^i+\ii y^i$ with $i=1,\dots, n$. 
By the coordinate transformation (\ref{coordinate-transform}), 
the bases are transformed by 
\begin{equation*}
 \fpartial{x_{(\mu)}} = \left( \varphi_{\lambda\mu}^t\right)^{-1} \fpartial{x_{(\lambda)}}, \qquad 
 \fpartial{x}:=(\fpartial{x^1},\dots,\fpartial{x^n})^t , 
\end{equation*}
and hence the corresponding coordinates are transformed by 
\begin{equation*}
 y_{(\mu)}=\varphi_{\lambda\mu} y_{(\lambda)},\qquad y:=(y^1,\dots, y^n)^t
\end{equation*}
so that the combination $\sum_i y^i\fpartial{x^i}$ is independent of 
the coordinate systems. 
This shows that 
the transition functions for the manifold $TB$ are given by 
\begin{equation*}
 \bp x_{(\mu)}\\ y_{(\mu)}\ep
 = \bp \varphi_{\lambda\mu} & 0 \\ 0 & \varphi_{\lambda\mu}\ep
 \bp x_{(\lambda)}\\ y_{(\lambda)}\ep
 + \bp \psi_{\lambda\mu} \\ 0 \ep, 
\end{equation*}
and hence the complex coordinate systems are transformed holomorphically: 
\begin{equation*}
 z_{(\mu)} = \varphi_{\lambda\mu} z_{(\lambda)} + \psi_{\lambda\mu} . 
\end{equation*}

On the other hand, for any smooth manifold $B$, the 
cotangent bundle $T^*B$ has a (canonical) symplectic form $\omega_{T^*B}$. 
For each $U_\lambda=U$, 
when we denote the coordinates of $T^*B|_U\simeq U\times\R^n$ 
by $(x^1,\dots,x^n;y_1,\dots,y_n)$, 
$\omega_{T^*B}$ is given by 
\begin{equation*}
 \omega_{T^*B}:= -d(\sum_{i=1}^n y_idx^i)= \sum_{i=1}^n dx^i\wedge dy_i . 
\end{equation*}
This is actually defined globally since 
the coordinate transformations on $T^*B$ are induced from the 
coordinate transformations of $\{U_\lambda\}_{\lambda\in\Lambda}$. 
Actually, one has 
\begin{equation*}
  dx_{(\lambda)}=\varphi_{\lambda\mu} dx_{(\mu)} 
\end{equation*}
and the corresponding coordinates are transformed by 
\begin{equation}\label{y_*transf}
 \check{y}_{(\lambda)}= \left(\varphi_{\lambda\mu}^t\right)^{-1} \check{y}_{(\mu)},\qquad 
\check{y}:=(y_1,\dots,y_n)^t
\end{equation}
so that the combination $\sum_{i=1}^n y_idx^i\in T^*B$ is independent of 
the coordinates. From this, it follows that the symplectic form 
$\omega_{T^*B}=-d(\sum_{i=1}^n y_idx^i)$ is defined globally.

By choosing a metric $g$ on a smooth manifold $B$, 
one obtains a bundle isomorphism between $TB$ and $T^*B$. 
For each $b\in B$, this isomorphism $TB\to T^*B$ is defined by 
$\xi\mapsto g(\xi, - )$ for $\xi\in T_bB$. 
This actually defines a bundle isomorphism since $g$ is nondegenerate 
at each point $b\in B$. 
This bundle isomorphism also induces a diffeomorphism 
from $TB$ to $T^*B$. 
In this sense, 
hereafter we sometimes identify $TB$ and $T^*B$. 
By this identification, $y^i$ and $y_i$ is related by 
\begin{equation*}
 y_i=\sum_{j=1}^n g_{ij}y^j , \qquad g_{ij}:=g\left(\fpartial{x^i},\fpartial{x^j}\right) . 
\end{equation*}

When an affine manifold $B$ is equipped with 
a metric $g$ which is expressed locally as 
\begin{equation*}
 g_{ij}
= \frac{\partial^2\phi}{\partial x^i\partial x^j}  
\end{equation*}
for some local smooth function $\phi$, 
then $(B,g)$ is called a {\em Hessian manifold}. 
When $B$ is a Hessian manifold, then $TB\simeq T^*B$ is equipped with 
the structure of K\"ahler manifold as we explain below. 
In this sense, 
a Hessian manifold is also called an affine K\"ahler manifold. 

First, when $B$ is affine, then $TB$ 
is already equipped with the complex structure $J_{TB}$. 
We fix a metric $g$ and set a two-form $\omega_{TB}$ on $TB$ as 
\begin{equation*}
 \omega_{TB}:=\sum_{i,j=1}^n g_{ij} dx^i\wedge  dy^j . 
\end{equation*}
This $\omega_{TB}$ is nondegenerate since $g$ is nondegenerate. 
Furthermore,
$\omega_{TB}$ is closed 
iff $(B,g)$ is Hessian, where $\omega_{TB}$ coincides with the 
pullback of $\omega_{T^*B}$ by the diffeomorphism $TB\to T^*B$. 
Thus, a Hessian manifold $(B,g)$ 
is equipped with the complex structure $J_{TB}$ and 
the symplectic structure $\omega_{TB}$. 
A metric $g_{TB}$ on $TB$ is then given by 
\begin{equation*}
 g_{TB}(X,Y):=\omega_{TB}(X,J_{TB}(Y)) 
\end{equation*}
for $X,Y\in\Gamma(T(TB))$. 
This is locally expressed as 
\begin{equation*}
 g_{TB} = \sum_{i,j=1}^n (g_{ij}dx^idx^j + g_{ij}dy^idy^j) . 
\end{equation*}
This shows that $g_{TB}$ is positive definite. 
To summarize, for a Hessian manifold $(B,g)$, $(TB, J_{TB},\omega_{TB})$ 
forms a K\"ahler manifold, where $g_{TB}$ is the K\"ahler metric. 

In order to define a K\"ahler structure on $T^*B$, 
we employ the dual affine local coordinate system 
\[
 \check{x}=(x_1,\dots, x_n)^t 
\]
of $x=(x^1,\dots, x^n)$ on $B$, that is, 
the coordinate system $\check{x}$ satisfying 
\begin{equation*}
 dx_i=\sum_{j=1}^n  g_{ij} dx^j . 
\end{equation*}
Such an $\check{x}$ actually exists since $(B,g)$ is Hessian; 
we may set  $x_i(x):=(\partial\phi/\partial x^i)(x)$. 
We thus obtain the dual coordinate system 
$\check{x}^{(\lambda)}:=(x_1^{(\lambda)},\dots,x_n^{(\lambda)})^t$  
for each $\lambda$. 
The dual coordinates 
then define another affine structure on $B$. 
Actually, 
the local description of the metric is changed by 
\begin{equation*}
 g_{(\lambda)}=\{(g_{(\lambda)})_{ij}\}_{i,j=1,\dots,n}
 =\left(\varphi_{\lambda\mu}^t\right)^{-1} g_{(\mu)}\varphi_{\lambda\mu}^{-1}, 
\end{equation*}
so one has $d\check{x}^{(\lambda)}=
\left(\varphi_{\lambda\mu}^t\right)^{-1}d\check{x}^{(\mu)}$ 
and then 
\begin{equation}\label{x_*transf}
 \check{x}^{(\lambda)} =
 \left(\varphi_{\lambda\mu}^t\right)^{-1} 
 \check{x}^{(\mu)} + \check{\psi}_{\lambda\mu}
\end{equation}
for some $\check{\psi}_{\lambda\mu}\in\R^n$. 
Thus, the combinations $z_i:=x_i+\ii y_i$, $i=1,\dots, n$, 
form a complex coordinate system on $T^*B$, and 
$T^*B$ forms a complex manifold. 
Actually, by eq.(\ref{y_*transf}) and (\ref{x_*transf}), 
one has the holomorphic coordinate transformation 
\begin{equation*}
 \check{z}^{(\mu)} 
 = \left(\varphi_{\lambda\mu}^t\right)^{-1} 
 \check{z}^{(\lambda)} + \check{\psi}_{\lambda\mu},\qquad 
\check{z}:=(z_1,\dots,z_n)^t . 
\end{equation*}

Using this dual coordinates, 
the symplectic form $\omega_{T^*B}$ is expressed locally as 
\begin{equation*}
 \omega_{T^*B}= \sum_{i,j=1}^n g^{ij} dx_i\wedge dy_j , 
\end{equation*}
where $g^{ij}$ is the $(i,j)$ element of the inverse matrix of 
$\{g_{ij}\}$. 
Then, we set a metric on $T^*B$ by 
\begin{equation*}
 g_{T^*B}(X,Y):=\omega_{T^*B}(X,J_{T^*B}(Y)) 
\end{equation*}
for $X,Y\in\Gamma(T(T^*B))$, 
which is locally expressed as 
\begin{equation*}
 g_{T^*B} = \sum_{i,j=1}^n (g^{ij}dx_idx_j + g^{ij}dy_idy_j) . 
\end{equation*}
These structures define a K\"ahler structure on $T^*B$.

For a tropical Hessian manifold $B$, 
we consider two $T^n$-fibrations over $B$ 
obtained by a quotient $M$ of $TB$ and 
a quotient $\check{M}$ of $T^*B$ 
by fiberwise $\Z^n$ action as follows. 

For $TB$, we locally consider $TB|_U$ and define 
a $\Z^n$-action generated by 
$y^i\mapsto y^i+ 2\pi$ for each $i=1,\dots, n$. 
For $T^*B$, we again locally consider $T^*B|_U$ and define 
a $\Z^n$-action generated by 
$y_i\mapsto y_i+ 2\pi$ for each $i=1,\dots, n$. 
Both $\Z^n$-actions are well-defined globally since 
$B$ is tropical affine, i.e., 
the transition functions of $n$-dimensional vector bundles $TB$ and 
$T^*B$ belong to $GL(n;\Z)$. 

Then, $M:=TB/\Z^n$ is a K\"ahler manifold whose 
symplectic structure $\omega_M$ and 
complex structure $J_M$ are those naturally induced from 
$\omega_{TB}$ and $J_{TB}$ on $TB$. 
Similarly, $\check{M}:=T^*B/\Z^n$ is a K\"ahler manifold 
whose symplectic structure $\omega_{\check{M}}$ and 
complex structure $J_{\check{M}}$ 
are those induced from $\omega_{T^*B}$ and $J_{T^*B}$, 
respectively.
In particular, $z=x+\ii y$ and $\check{z}=\check{x}+\ii\check{y}$ 
turn out to be local complex coordinates of the complex manifolds $M$ and $\check{M}$, respectively. 
The fibrations $\pi:M\to B$ and $\pi\check{}:\check{M}\to B$ are often called 
{\em semi-flat torus fibrations} or {\em $T^n$-invariant manifolds}. 
See \cite{LYZ,leung05} and also \cite{fukaya:asymptotic}. 
Since $M$ and $\check{M}$ are dual to each other, 
we can construct them in the opposite way. That is,
if we consider the coordinate systems 
${\check{x}^{(\lambda)}}$ for
$B$, then the tangent bundle over $B$ is $T^*B$ above, and
the cotangent bundle is $TB$. 
Following \cite{LYZ,leung05}, we treat $M$ as a symplectic manifold
and $\check{M}$ as a complex manifold and
discuss the homological mirror symmetry.

 \subsection{Toric manifolds and $T^n$-invariant manifolds}
\label{ssec:tt}
 
The set-up in the previous subsection is originally applied to 
the mirror symmetry of compact Calabi-Yau manifolds $M,\check{M}$. 
We would like to extend this set-up 
to the case $\check{M}$ is the complement of the toric divisors
of a smooth compact toric manifold $X$. 
The complement $\check{M}$ is actually a trivial torus fibration
$\pi\check{}:\check{M}\to B$ where the base $B$ is identified with 
the interior of the dual polytope $P$ of $X$. 

What may be more interesting is that
$B$ is actually tropical affine in this situation.
Of course, since $B$ is a contractible open set, $B=\Int(P)$ has
an open covering by itself, which means that $B$ is tropical affine. 
However, what we meant is something stronger in the following sense.
A smooth compact toric manifold $X$ has a natural open covering 
$\{\ti{\cU}_\lambda\}_{\lambda\in\Lambda}$, where each 
$\ti{\cU}_\lambda$ is associated to each cone in the fan of $X$ of maximal dimension, 
which induces 
the open covering 
$\{\ti{U}_\lambda:=\pi\check{}(\ti{\cU}_\lambda)\}_{\lambda\in\Lambda}$ of $P$ 
(where the origin of each $\ti{\cU}_\lambda$ corresponds by $\pi\check{}$ to each vertex of $P$). 
Then we see that the coordinate transformations
are tropical affine 
(though $\ti{U}_\lambda\cap B=B$ for any $\lambda$.)
This seems important
since we need to include some information from the boundary
$\partial(B)$ of $B$ when we discuss homological mirror symmetry of $X$
and its mirror dual.

Let us see the above construction explicitly for
$X=\C P^n$. 
For
\[
 \C P^n = \{[t_0:\cdots: t_n]\} , 
\]
the natural open covering is $\{\ti{\cU}_\lambda\}_{\lambda=0,1,...,n}$ where
\[
 \ti{\cU}_\lambda=\{[t_0:\cdots: t_n]\ |\ t_\lambda\ne 0\} . 
\]
The corresponding local coordinates are $(w^{(\lambda)}_1,\dots, w^{(\lambda)}_n)$ where 
\begin{equation}\label{w}
w^{(\lambda)}_1=t_0/t_\lambda,\ ... ,\ w^{(\lambda)}_{\lambda}=t_{\lambda-1}/t_\lambda,\ 
w^{(\lambda)}_{\lambda+1}=t_{\lambda+1}/t_\lambda,\ ...,\
w^{(\lambda)}_n=t_n/t_\lambda . 
\end{equation}
We identify $\check{M}$ with the complement of the toric divisors
of $\C P^n$: 
\[
\check{M} = \{ [t_0:\cdots : t_n]\ |\ t_0\cdot t_1\cdots t_n\ne 0\} ,
\]
where $\pi\check{}: \check{M}\to B$ is given by
\[
 \pi\check{}\, ([t_0:\cdots : t_n]):=[|t_0|:\cdots : |t_n|] 
\]
(though we express this in a different coordinate system below). 
So we have $\cU_\lambda:=\ti{\cU}_\lambda\cap\check{M}= \check{M}$
for any $\lambda$. 
We further denote $U_\lambda:=\pi\check{}\, (\cU_\lambda)$. 
For each $\cU_\lambda$, we express
\[
 w^{(\lambda)}_i=: e^{z^{(\lambda)}_i}=e^{x^{(\lambda)}_i+\ii y^{(\lambda)}_i} .
\]
Since the coordinate transformation between
$\check{z}^{(\lambda)}$ and $\check{z}^{(\mu)}$
is tropical affine by (\ref{w}), so is 
the coordinate transformation between
$\check{x}^{(\lambda)}$ and $\check{x}^{(\mu)}$. 

Hereafter we consider $U:=U_0$ 
(since $U_0=U_1=\cdots =U_n=B$) 
and drop the upper index ${}^{(0)}$; for instance
$w_i^{(0)}=:w_i$ and $x^{(0)}_i=:x_i$. 
The Fubini-Study K\"ahler form is then expressed in $\check{M}=(\pi\check{}\, )^{-1}(U)$ as
\[
\omega_{\check{M}}=-2\ii d\left(
 \frac{\bar{w}_1dw_1+\cdots +\bar{w}_ndw_n}{1+\bar{w}_1w_1+\cdots +\bar{w}_nw_n}\right) . 
\]
When we express this as 
$\omega_{\check{M}}=\sum_{i,j=1}^n g^{ij} dx_i \wedge dy_j$, we have 
\[
 \begin{split}
  & g^{ij}=\frac{\partial^2\check{\phi}}{\partial x_i \partial x_j},  \\
   & \check{\phi}=\log(1+e^{2x_1}+\cdots + e^{2x_n}) .
 \end{split}
\]
Thus, $B$ is a Hessian manifold.  
The dual coordinates $(x^1,\dots, x^n)$ is obtained by 
\[
dx^i= \sum_{j=1}^n \frac{\partial^2\check{\phi}}{\partial x_i \partial x_j} dx^j
    = d\left(\fpart{\check{\phi}}{x_i}\right), 
\]
so
\begin{equation}\label{dual-x}
 x^i=\fpart{\check{\phi}}{x_i} = \frac{2e^{2x_i}}{1+e^{2x_1}+\cdots + e^{2x_n}} .
\end{equation}
By this $(x^1,\dots, x^n)$, $B$ is expressed as
\[
 B=\{(x^1,\dots, x^n)|\ x^1>0,\dots, x^n>0, x^1+\cdots + x^n<2\} . 
\]
In particular, in this coordinate system, $\pi\check{}: \check{M}\to B$ is expressed as 
\[
 \begin{split}
 \pi\check{}(\check{x},\check{y}) & = (x^1(\check{x}),\dots, x^n(\check{x})) \\
 & \left( = \left( \frac{2\bar{w}_1w_1}{1+\bar{w}_1w_1+\cdots +\bar{w}_nw_n},\dots,\frac{2\bar{w}_nw_n}{1+\bar{w}_1w_1+\cdots +\bar{w}_nw_n}\right) \right). 
 \end{split}
\]
In this way, $\pi\check{}$ is regarded as the restriction to $\check{M}\subset X$ 
of the moment map 
$X\to\R^n$ for the $T^n$ action on $X = \C P^n$. 

Note also that we can regard $M$, the dual torus fibration of $\check{M}$,
as an open complex submanifold of $(\C^\times)^n$, 
where $(e^{-(x^1+\ii y^1)}, \dots, e^{-(x^n+\ii y^n)})$ is the coordinate system
of $(\C^\times)^n$. 
However, the symplectic form $\omega_M=\sum_{i,j=1}^n g_{ij}dx^i\wedge dy^j$ diverges at the boundary $\partial(M)=\pi^{-1}(\partial(B))$.

 \section{Lagrangian submanifolds and holomorphic vector bundles}
\label{sec:Lag-Vect}

In the first two subsections, we first recall
the construction of line bundles on $\check{M}$
associated to Lagrangian sections of $M\to B$ 
discussed in \cite{LYZ,leung05}. 
Then, in subsection \ref{ssec:cpn-Lag}, we apply this construction to
the case $\check{M}$ is the complement of the
toric divisors of $\C P^n$. 

This construction gives an objectwise correspondence of the corresponding homological mirror symmetry. 
More generally, on $M$ we can consider Lagrangian sections equipped with local systems 
as objects of the Fukaya category $Fuk(M)$. 
However, we do not discuss this generalized set-up since we need only Lagrangian sections equipped with trivial local system for our purpose. See also Remark \ref{rem:local-system} at the end of subsection \ref{ssec:Vect}.

 \subsection{Lagrangian submanifolds in $M$}
\label{ssec:Lag}

We fix a tropical affine open covering $\{U_\lambda\}_{\lambda\in\Lambda}$. 
Let $\sbar: B\to M$ be a section of $M\to B$. 
Locally, we may regard $\sbar$ as a section of $TB\simeq T^*B$. 
Then, $\sbar$ is locally described by a collection of functions as 
\begin{equation*}
 y^i_{(\lambda)}=s^i_{(\lambda)}(x) 
\end{equation*}
on each $U_\lambda$. 

On $U_\lambda\cap U_\mu$, these local expressions 
are related to each other by 
\begin{equation}\label{s-transform}
 s_{(\mu)}(x) = s_{(\lambda)}(x) + 2\pi I_{\lambda\mu} 
\end{equation}
for some $I_{\lambda\mu}\in\Z^n$. 
Here, 
$x$ may be identified with either $x_{(\lambda)}$ 
or $x_{(\mu)}$. Also, $s_{(\lambda)}(x)$ and $s_{(\mu)}(x)$ are 
expressed by the common coordinates $y_{(\lambda)}$ or $y_{(\mu)}$. 
This transformation rule automatically satisfies the 
cocycle condition 
\begin{equation}\label{cocycle-cd}
I_{\lambda\mu}+I_{\mu\nu}+I_{\nu\lambda}=0
\end{equation}
for $U_\lambda\cap U_\mu\cap U_\nu\ne\emptyset$. 
We denote by $s$ such a collection 
$\{s_{(\lambda)}: U_\lambda\to TB|_{U_\lambda}\}_{\lambda\in\Lambda}$ 
which is equipped with the transformation rule (\ref{s-transform}) 
satisfying the cocycle condition (\ref{cocycle-cd}). 

Now we discuss when the graph of $\sbar$ forms a Lagrangian submanifold 
in $M$. 
By definition, an $n$-dimensional submanifold $L$ in a 
$2n$-dimensional symplectic manifold $(M,\omega_M)$ is 
{\em Lagrangian} iff $\omega_M|_L=0$. 
This is a local condition. 
Thus, in order to discuss whether 
the graph of a section 
$\sbar: B\to M$ is Lagrangian or not, 
we may check the condition locally and in particular in $T^*B$. 

It is known (as shown easily by taking the basis) 
that the graph of 
$\sum_{i=1}^n y_idx^i$ with local functions $y_i$ is Lagrangian in $T^*B$ 
iff there exists a local function $f$ such that
$\sum_{i=1}^n y_idx^i=df$. 
Now, a section $\sbar :B\to M$ is locally regarded as a section of 
$T^*B$ by setting 
$y_i = \sum_{j=1}^n g_{ij}y^j = \sum_{j=1}^n g_{ij} s^j$, from which one has 
\begin{equation*}
 \sum_{i=1}^n y_idx^i=\sum_{i=1}^n (\sum_{j=1}^n g_{ij}s^j) dx^i=\sum_{j=1}^n s^jdx_j. 
\end{equation*}
Thus, the graph of the section $\sbar:B\to M$ is Lagrangian 
iff there exists a local function $f$ such that 
$\sum_{j=1}^n s^j dx_j=df$. 

Note that $y=s(x)$ defines a special Lagrangian submanifold
if $s$ is affine with respect to $x^i$. 
(Thus, the zero section of $M\to B$ 
is a special Lagrangian submanifold. )

The gradient vector field is of the form:
\begin{equation}\label{grad-f}
 \grad(f):= \sum_{i,j}\fpart{f}{x^j}g^{ji}\fpartial{x^i}=
 \sum_i \fpart{f}{x_i}\fpartial{x^i} . 
\end{equation}

 \subsection{Holomorphic vector bundles on $\check{M}$}
\label{ssec:Vect}

Consider a section $\sbar: B\to M$ and express it 
as a collection $s=\{s_{(\lambda)}\}_{\lambda\in\Lambda}$ of local functions. 
We define a line bundle $V$ with a $U(1)$-connection 
on the mirror manifold $\check{M}$ 
associated to $s$. 
We set the covariant derivative locally as
\footnote{We switch the sign of the connection one form compared to
that in \cite{hk:fukayadeform} so that 
the mirror correspondence of objects fits with the one
in homological mirror symmetry of tori as in \cite{hk:revisited} 
and references therein. 
We also include $2\pi$ in various places which are missing in \cite{hk:fukayadeform}. 
}
\begin{equation}\label{conn}
  D:=d-\frac{\ii}{2\pi} \sum_{i=1}^n s^i(x) dy_i , 
\end{equation}
whose curvature is 
\begin{equation*}
 D^2= \frac{\ii}{2\pi}\sum_{i,j=1}^n \fpart{s^i}{x_j} dx_j\wedge dy_i. 
\end{equation*}
The $(0,2)$-part vanishes iff 
the matrix $\fpart{s^i}{x_j}$ is symmetric, which is the case 
when there exists a function $f$ locally such that 
$df = \sum_{i=1}^n s^idx_i$. 
Thus, the condition that 
$D$ defines a holomorphic line bundle on $\check{M}$ is 
equivalent to that the graph of $\sbar$ is Lagrangian in $M$.

This covariant derivative $D$ is in fact defined globally. 
Suppose that $D$ is given locally on each $\check{M}|_{U_\lambda}$ 
of the $T^n$-fibration $\check{M}\to B$ 
with a fixed tropical affine open covering $\{U_\lambda\}_{\lambda\in\Lambda}$. 
Namely, we continue to employ $\{U_\lambda\}_{\lambda\in\Lambda}$ 
for local trivializations of the line bundle associated to a section 
$\sbar: B\to M$. 
The transition functions for $(V,D)$ are 
defined as follows. 
Recall that the section $\sbar: B\to M$ is expressed locally 
as 
\begin{equation*}
 y^i_{(\lambda)}=s^i_{(\lambda)}(x) 
\end{equation*}
on each $U_\lambda$, where, on $U_\lambda\cap U_\mu$, the local expression 
is related to each other by 
\begin{equation*}
 s_{(\mu)}(x) = s_{(\lambda)}(x) + 2\pi I_{\lambda\mu} 
\end{equation*}
for some $I_{\lambda\mu}\in\Z^n$ (see eq.(\ref{s-transform})). 
Correspondingly, the transition function for the line bundle $V$ with 
the connection $D$ is given by 
\begin{equation*}
 \psi_{(\mu)} = e^{\ii I_{\lambda\mu}\cdot\check{y}}\psi_{(\lambda)}  
\end{equation*}
for local expressions $\psi_{(\lambda)}$, $\psi_{(\mu)}$ of 
a smooth section $\psi$ of $V$, where 
$I_{\lambda\mu}\cdot\check{y}:=\sum_{j=1}^n i_{j}y_j$ for 
$I_{\lambda\mu}=(i_1,\dots,i_n)$. 
We see the compatibility 
\begin{equation*}
 (D\psi_{(\lambda)})_{(\mu)} = D(\psi_{(\mu)})  
\end{equation*}
holds true since 
the left hand side turns out to be 
\begin{equation*}
 \begin{split}
 & e^{\ii I_{\lambda\mu}\cdot\check{y}}
 \left(\left(d-\frac{\ii}{2\pi} s_{(\lambda)}(x) \cdot dy\right)
 e^{-\ii I_{\lambda\mu}\cdot\check{y}}\psi_{(\mu)}\right) \\
 & = e^{\ii I_{\lambda\mu}\cdot\check{y}} e^{-\ii I_{\lambda\mu}\cdot\check{y}} 
 \left(\left(d-\frac{\ii}{2\pi} (s_{(\lambda)}(x)+2\pi I_{\lambda\mu})\cdot dy\right)\psi_{(\mu)} \right) \\
 & = \left(d-\frac{\ii}{2\pi} s_{(\mu)}(x)\cdot dy\right)\psi_{(\mu)} . 
 \end{split}
\end{equation*}
Since $(V,D)$ is locally-trivialized 
by $\{\check{M}|_{U_\lambda}\}_{\lambda\in\Lambda}$, 
for each $x\in B$, 
$\psi(x, \cdot )$ gives a smooth function on the fiber $T^n$. 
Thus, on each $U_\lambda$, $\psi(x,y)$ can be Fourier-expanded as 
\begin{equation*}
 \psi(x,y)|_{U_\lambda}=\sum_{I\in\Z^n} \psi_{\lambda,I}(x)e^{\ii I\cdot\check{y}},  
\end{equation*}
where $I\cdot\check{y}:=\sum_{j=1}^n i_j y_j$ for $I=(i_1,\dots,i_n)$. 
Note that each coefficient $\psi_{\lambda,I}$ is a smooth function on $U_\lambda$. 
In this expression, the transition function acts to each 
$\psi_{\lambda,I}$ as 
\begin{equation*}
 \begin{split}
 \sum_{I\in\Z^n}\psi_{\mu,I} e^{\ii I\cdot\check{y}} & = 
 e^{\ii I_{\lambda\mu}\cdot\check{y}} \sum_{I\in\Z^n} \psi_{\lambda,I} e^{\ii I\cdot\check{y}} \\
 &= \sum_{I\in\Z^n} \psi_{\lambda,I} e^{\ii (I+I_{\lambda\mu})\cdot\check{y}} \\
 &= \sum_{I\in\Z^n} \psi_{\lambda,I-I_{\lambda\mu}} e^{\ii I\cdot\check{y}} \\
 \end{split}
\end{equation*}
and hence $\psi_{\mu,I}=\psi_{\lambda,I-I_{\lambda\mu}}$. 
\begin{rem}\label{rem:local-system}
We can also associate a line bundle to a Lagrangian section {\em equipped with a local system}, 
where the holonomy turns out to be included as the real coefficient of 
each $dy_i$ in the covariant derivation (\ref{conn}). 
We do not discuss this generalized set-up 
since the real coefficients are trivial for any line bundle we need in this paper 
and hence the corresponding object is a Lagrangian section with the trivial local system.  
\end{rem}

\subsection{Holomorphic line bundles on {$\C P^n$} and the corresponding Lagrangians}
\label{ssec:cpn-Lag}

In the previous subsections, 
we assign a line bundle on $\check{M}$ 
to each Lagrangian section in $M\to B$.
In this subsection, we start from a line bundle on $\C P^n$. 
We identify $\check{M}$
with the complement of the toric divisors of $\C P^n$,
and restrict the line bundle to $\check{M}$. 
We will see that, by twisting it with an appropriate isomorphism, 
the result actually comes from a Lagrangian section in $M\to B$. 
In this way, we construct a Lagrangian section in $M\to B$ 
corresponding to $\cO(a)$ on $\C P^n$ for any $a\in\Z$. 

We continue the convention in subsection \ref{ssec:tt}. 
The complement of the toric divisors of $\C P^n$
is
\[
 \check{M}=\{[t_0:\cdots: t_n]\ |\ t_0\cdot t_1\cdots t_n\ne 0\}, 
\]
where 
\[
 e^{x_i+\ii y_i} = w_i = t_i/t_0 . 
\]
A connection one-form of $\cO(a)$ is given by the one
which is expressed locally on $\check{M}$ as  
\begin{equation}\label{conn-cpn}
 \begin{split}
A_a & = -a
\frac{\bar{w}_1dw_1+\cdots + \bar{w}_ndw_n}{1+\bar{w}_1w_1+\cdots +\bar{w}_nw_n} \\
 & = -a
\frac{e^{2x_1}(dx_1+\ii dy_1)+\cdots + e^{2x_n}(dx_n+\ii dy_n)}{1+e^{2x_1}+\cdots +e^{2x_n}} .
 \end{split}
\end{equation}
We twist this by 
\[
 \Psi_a:= (1+e^{2x_1}+\cdots + e^{2x_n})^{a/2}, 
\]
and then obtain 
\begin{equation}\label{Psi_a}
 \Psi_a^{-1}(d+A_a)\Psi_a= d-\ii a
 \frac{e^{2x_1}dy_1+\cdots + e^{2x_n}dy_n}{1+e^{2x_1}+\cdots +e^{2x_n}} .
\end{equation}
By the previous subsections, this is the line bundle on
$\check{M}$ which corresponds to
the Lagrangian section $L_a$ in $M\to B$ expressed as 
\[
 \bp y^1 \\ \vdots \\ y^n \ep 
 = \bp s_a^1 \\ \vdots \\ s_a^n \ep  
 = 2\pi a \bp \frac{e^{2x_1}}{1+e^{2x_1}+\cdots + e^{2x_n}} \\ \vdots \\
  \frac{e^{2x_n}}{1+ e^{2x_1}+\cdots + e^{2x_n}}\ep 
 =2\pi \frac{a}{2} \bp x^1 \\ \vdots \\ x^n \ep  
\]
where $x^i>0$ for $i=1,\dots, n$ and $x^1+\cdots + x^n< 2$. 

We see that the local function 
\begin{equation*}
  \begin{split}
  f_a& =2\pi \frac{a}{2}(\log (1+e^{2x_1}+\cdots + e^{2x_n}) -\log 2)  \\
   & =-2\pi \frac{a}{2}\log (2-x^1-x^2-\cdots -x^n)
 \end{split}
\end{equation*}
satisfies $df_a= \sum_{i=1}^n s^i_a dx_i$. 
The corresponding gradient vector field is
\begin{equation*}
 \grad(f_a) = \sum_{i=1}^n \fpart{f_a}{x_i}\fpartial{x^i} 
 = 2\pi \frac{a}{2} \left(x^1\fpartial{x^1} +\cdots + x^n\fpartial{x^n}\right) 
\end{equation*}
by (\ref{grad-f}). 
\begin{rem}
This Lagrangian section $L_a$ is a special Lagrangian since
it is expressed locally as the graph of linear functions $y^i(x)$
of $(x^1,\dots, x^n)$. 

Furthermore, we see that $L_a$ includes a critical point of 
the corresponding Landau-Ginzburg potential. 
In fact, the Landau-Ginzburg potential $W$ is 

\[
 W(w^1,w^2,\dots,w^n) := w^1+\cdots + w^n + \frac{e^{-2}}{w^1w^2\cdots w^n} , 
\]
where $w^i:=e^{-(x^i+\ii y^i)}$. 
The critical points are given by
\[
 (w^1,\dots, w^n)= \left(
 e^{-\frac{2}{n+1}}\zeta^a,\dots,
 e^{-\frac{2}{n+1}}\zeta^a\right)=:c_a ,\qquad a=0, 1,\dots, n
\]
where $\zeta=(1)^{\frac{1}{n+1}}$ is the $(n+1)$-th root of unity. 
Thus, we see that
each critical point $c_a\in (\C^\times)^n$ is included in 
$L_a$. 
\end{rem}

\section{Homological mirror symmetry set-up}
\label{sec:Cat}
 
In this section,
we first recall DG-categories $\cF(M)$ and $\cV(\check{M})$ 
associated to
$M$ and $\check{M}$ respectively, following \cite{hk:fukayadeform}.
Kontsevich-Soibelman's approach for the homological mirror symmetry \cite{KoSo:torus}
introduces an intermediate category $\Mo(B)$ and the existence of
an $A_\infty$-equivalence
\[
 Fuk(M) \simeq \Mo(B) \overset{\sim}{\to} \cF(M) .
\]
This $\Mo(B)$ is called the weighted Fukaya-Oh category or
the category of weighted Morse homotopy on $B$. 
In subsection \ref{ssec:Mo}, we propose a modification $\Mo(P)$ of $\Mo(B)$ 
where $P=\bar{B}$ is the dual polytope of a smooth compact toric manifold
$X$.

 \subsection{DG-category $\cV$ associated to $\check{M}$}
\label{ssec:DG_cM}

We define a DG-category $\cV=\cV(\check{M})$
of holomorphic line bundles over $\check{M}$ as follows. 
The objects are holomorphic line bundles $V$ with $U(1)$-connections $D$ 
associated to lifts $s$ of sections 
as we defined in subsection \ref{ssec:Vect}. 
We sometimes label these objects as $s$ instead of $(V,D)$. 
For any two objects $s_a=(V_a,D_a), s_b=(V_b,D_b)\in\cV$, 
the space $\cV(s_a,s_b)$ of morphisms 
is defined by 
\begin{equation*}
 \cV(s_a,s_b):= 
 \Gamma(V_a,V_b)\otimes_{C^\infty(\check{M})}\Omega^{0,*}(\check{M}), 
\end{equation*}
where 
$\Omega^{0, *}(\check{M})$ is the space of 
anti-holomorphic differential forms 
and $\Gamma(V_a,V_b)$ is the space of homomorphisms 
from $V_a$ to $V_b$.
\footnote{Here we again make a minor change 
of the formulation of the DG category compared to \cite{hk:fukayadeform}
due to the change of sign in (\ref{conn}). 
}
The space $\cV(s_a,s_b)$ is a $\Z$-graded vector space, 
where the grading is defined 
as the degree of the anti-holomorphic differential forms. 
The degree $r$ part is denoted by $\cV^r(s_a,s_b)$. 
We define a linear map 
$d_{ab}:\cV^r(s_a,s_b)\to\cV^{r+1}(s_a,s_b)$ 
as follows. We decompose $D_a$ into its holomorphic part 
and anti-holomorphic part $D_a=D_a^{(1,0)}+D_a^{(0,1)}$, and 
set $2 D_a^{(0,1)}=:d_a$. Then, for $\psi\in\cV^r(s_a,s_b)$, 
we set 
\begin{equation*}
 d_{ab}(\psi):= d_b \psi - (-1)^r \psi d_a\in\cV^{r+1}(s_a,s_b) . 
\end{equation*}
Note that $d_{ab}^2=0$ since each $(V_a,D_a)$ is holomorphic, i.e.,
$(d_a)^2=0$. 

The product structure 
$m: \cV(s_a,s_b)\otimes \cV(s_b,s_c)\to \cV(s_a,s_c)$ 
is defined by the composition 
of homomorphisms of line bundles together with the 
wedge product for the anti-holomorphic differential forms. 
More precisely,
for $\psi_{ab}\in\cV^{r_{ab}}(s_a,s_b)$ and
$\psi_{bc}\in\cV^{r_{bc}}(s_b,s_c)$, we set 
\[
m(\psi_{ab}, \psi_{bc}):=(-1)^{r_{ab}r_{bc}}\psi_{bc}\wedge\psi_{ab} 
\ (=\psi_{ab}\wedge\psi_{bc}),
\]
where $\wedge$ denotes the operation consisting of the composition and the
wedge product. 
Then, we see that $\cV$ forms a DG-category. 
\footnote{
In \cite{hk:fukayadeform}, we construct a curved DG category $DG_{\check{M}}$
where the objects are not necessarily holomorphic. 
The relation is given by $\cV=DG_{\check{M}}(0)$. 
}

In order to construct another equivalent curved DG-category, 
we rewrite this DG-category $\cV$ more explicitly. 
For an element $\psi\in \cV^r(s_a,s_b)$, 
we Fourier-expand this locally as 
\begin{equation*}
\psi(\check{x},\check{y})=\sum_{I\in\Z^n} \psi_I(\check{x}) e^{ \ii I\cdot \check{y}}, 
\end{equation*}
where $\psi_I$ is a smooth anti-holomorphic differential form 
of degree $r$. Namely, it is expressed as 
\begin{equation*}
 \psi_I= \sum_{i_1,\dots,i_r}\psi_{I;i_1\cdots i_r} d\zb_{i_1}\wedge\cdots\wedge d\zb_{i_r}
\end{equation*}
with smooth functions $\psi_{I;i_1\cdots i_r}$. 
Let us express the transformation rule for $s_a$ as 
\begin{equation}
 (s_a)_{(\mu)}=(s_a)_{(\lambda)}+2\pi I_a \label{transf:s}
\end{equation}
with $I_a=I_{a;\lambda\mu}\in\Z^n$. 
The transition function is then given by 
$\psi_{(\mu)}= e^{\ii (I_b-I_a)\cdot\check{y}}\psi_{(\lambda)}$ and hence 
\begin{equation*}
 \psi_{(\mu),I}= \psi_{(\lambda),I+I_a-I_b} .  
\end{equation*}

The differential $d_{ab}$ is expressed locally as follows. 
Since 
\begin{equation*}
 \begin{split}
 D_a & =d-\frac{\ii}{2\pi} \sum_{j=1}^n s_a^j(x) dy_j \\
     & = \sum_{j=1}^n
 \left( \fpartial{x_j}dx_j+\left(\fpartial{y_j}- \frac{\ii}{2\pi} s_a^j\right)dy_j\right) \\
 & = \ov{2}\sum_{j=1}^n\left(\fpartial{x_j}-\ii\left(\fpartial{y_j}-\frac{\ii}{2\pi} s_a^j\right)\right)dz_j
 + \ov{2}\sum_{j=1}^n\left(\fpartial{x_j}+\ii\left(\fpartial{y_j}-\frac{\ii}{2\pi} s_a^j\right)\right)d\zb_j, 
 \end{split}
\end{equation*}
one has 
\begin{equation*}
 d_{a}=2D_a^{(0,1)}=
 \sum_{j=1}^n\left(\fpartial{x_j}+ \frac{s_a^j}{2\pi}+\ii\fpartial{y_j}\right)d\zb_j 
\end{equation*}
and then 
\begin{equation}\label{d_ab-local}
 d_{ab}(\psi)= 2\bpartial(\psi) - \frac{1}{2\pi}\sum_{i=1}^n(s_a-s_b)^id\zb_i\wedge \psi .  
\end{equation}

\subsection{DG-category $\cF$ associated to $M$}
\label{ssec:DG_M}

We define a DG-category $\cF=\cF(M)$ consisting of Lagrangian sections in $M$
as follows
\footnote{This $\cF$ corresponds to $DG_M(0)$ in
\cite{hk:fukayadeform}.}. 
As we shall see, we construct it so that it is canonically isomorphic 
to the previous DG-category $\cV$. 
We fix a tropical affine open covering 
$\{U_\lambda\}_{\lambda\in\Lambda}$ of $B$. 

After the identification of $(V,D)$ with $s$ 
made in the previous subsection \ref{ssec:DG_cM}, 
the objects are the same as those in $\cV$, that is, 
lifts $s$ of sections of $M\to B$. 
Under this identification the object $s_a\in \cF$, 
satisfies the transformation rule (\ref{transf:s}) as above. 
For each $\lambda\in\Lambda$ and $I\in\Z^n$, 
let $\Omega_{\lambda,I}(s_a,s_b)$ be the space of 
complex valued smooth differential forms on $U_\lambda$. 
The space $\cF(s_a,s_b)$ is then the subspace of 
\begin{equation*}
 \prod_{\lambda\in\Lambda}\prod_{I\in\Z^n} \Omega_{\lambda,I}(s_a,s_b)
\end{equation*}
such that 
\begin{itemize}
 \item $\phi_{\lambda,I}\in\Omega_{\lambda,I}(s_a,s_b)$ satisfies 
\begin{equation*}
 \phi_{\mu,I}|_{U_\lambda\cap U_\mu}= \phi_{\lambda,I+I_a-I_b}|_{U_\lambda\cap U_\mu} 
\end{equation*}
for any $U_\lambda\cap U_\mu\ne\emptyset$ and 
 \item the sum $\sum_{I\in\Z^n} \phi_{\lambda,I} e^{\ii I\cdot\check{y}}$ 
converges as smooth differential forms on each $M|_{U_\lambda}$. 
\end{itemize}
The space $\cF(s_a,s_b)$ is a $\Z$-graded vector space, 
where the grading is defined 
as the degree of the differential forms. 
The degree $r$ part is denoted $\cF^r(s_a,s_b)$. 
We define a linear map 
$d_{ab}:\cF^r(s_a,s_b)\to\cF^{r+1}(s_a,s_b)$ 
which is expressed locally as 
\begin{equation*}
 d_{ab}(\phi_{\lambda,I}):=d(\phi_{\lambda,I}) 
 -\frac{1}{2\pi}\sum_{j=1}^n (s^j_a-s^j_b + 2\pi i_j)dx_j\wedge \phi_{\lambda,I}
\end{equation*}
for $\phi_{\lambda,I}\in\Omega_{\lambda,I}(s_a,s_b)$ with $I:=(i_1,\dots,i_n)\in\Z^n$, 
where $d$ is the exterior differential on $B$. 
We have $d_{ab}^2=0$. 

The composition of morphisms 
$m: \cF(s_a,s_b)\otimes \cF(s_b,s_c)\to \cF(s_a,s_c)$ 
is defined by 
\begin{equation*}
 m(\phi_{ab;\lambda,I},\phi_{bc;\lambda,J}):=\phi_{ab;\lambda,I}\wedge\phi_{bc;\lambda,J}\ 
 \in\Omega_{\lambda,I+J}(s_a,s_c)
\end{equation*}
for $\phi_{ab;\lambda,I}\in\Omega_{\lambda,I}(s_a,s_b)$ and 
$\phi_{bc;\lambda,J}\in\Omega_{\lambda,I}(s_b,s_c)$. 
These structures define a DG-category $\cF$. 
Note that this $\cF$ is believed to be $A_\infty$-equivalent 
to the corresponding full subcategory of the Fukaya 
category $Fuk(M)$. 
(Compare this $\cF$ with 
what is called 
the deRham model for the Fukaya category
in Kontsevich-Soibelman \cite{KoSo:torus}, 
in particular a construction in the Appendix (Section 9.2). ) 
In subsection \ref{ssec:compare}, we shall explain the outline of 
how to compare $\cF$ with the Fukaya category.

 \subsection{Equivalence between $\cF$ and $\cV$}
\label{ssec:equivalent}

The DG-category $\cF$ is canonically isomorphic to 
the DG-category $\cV$. 
In fact, we see that the objects in $\cF$ are the same 
as those in $\cV$. 
The spaces of morphisms 
in $\cF$ and in $\cV$ are also identified canonically as follows. 
For a morphism $\phi_{ab}=\{\phi_{ab;\lambda,I}\}\in\cF^r(s_a,s_b)$, 
each $\phi_{ab;\lambda,I}$ is expressed as 
\begin{equation*}
 \phi_{ab;\lambda,I}= \sum_{i_1,\dots,i_r} \phi_{ab;\lambda,I;i_1\cdots i_r} 
 dx_{i_1}\wedge\cdots\wedge dx_{i_r}. 
\end{equation*}
To this, we correspond an element in $\cV^r(s_a,s_b)$ 
which is locally given as 
\begin{equation*}
 \sum_{i_1,\dots,i_r} (\phi_{ab;\lambda,I;i_1\cdots i_r} e^{\ii I\cdot\check{y}})
 d\zb_{i_1}\wedge\cdots\wedge d\zb_{i_r} 
\end{equation*}
on $U_\lambda$. 
We denote this correspondence by $\f:\cF\to \cV$. 
It is easily seen that our construction guarantees the following fact. 
\begin{prop}[{\cite[Proposition 4.1.]{hk:fukayadeform}}]
The functor $\f:\cF\to\cV$ 
is a DG-isomorphism. 
\end{prop}

 \subsection{The DG-category $\cF$ and the Fukaya category $Fuk(M)$}
\label{ssec:compare}

The DG-category $\cF$ is expected to be $A_\infty$-equivalent to
the Fukaya category $Fuk(M)$ \cite{fukayaAinfty}
of Lagrangian sections. 
The idea discussed in \cite{KoSo:torus}
to relate them is 
to apply homological perturbation theory 
to the DG-category $\cF$ (as an $A_\infty$-category) 
in an appropriate way so that the induced $A_\infty$-category 
coincides with 
(the full subcategory of) the Fukaya category $Fuk(M)$. 
More precisely, what should be induced directly from $\cF$ 
is the category $\Mo(B)$ of weighted Morse homotopy or
the Fukaya-Oh category for the torus fibration $M\to B$ introduced 
in section 5.2 of \cite{KoSo:torus}. 
Here, the Fukaya-Oh category means the $A_\infty$-category 
of Morse homotopy on $B$ introduced in \cite{fukayaAinfty}. 
It is shown in \cite{FO} that the Fukaya-Oh category is 
equivalent to the Fukaya category $Fuk(T^*B)$ 
consisting of the corresponding objects. 
The Fukaya-Oh category for the torus fibration $M\to B$ 
is a generalization of the Fukaya-Oh category on $B$ 
so that it corresponds to 
the Fukaya category $Fuk(M)$ instead of $Fuk(T^*B)$. 
Thus, a natural way of obtaining the $A_\infty$-equivalence $\cF\simeq Fuk(M)$
is to interpolate the category $\Mo(B)$ so that
\[
 Fuk(M) \simeq \Mo(B) \overset{\sim}{\to} \cF ,  
\]
where an $A_\infty$-equivalence $\Mo(B)\to\cF$ is expected to be obtained by
the homological perturbation theory 
and then $\Mo(B)$ is identified with $Fuk(M)$ by a linear $A_\infty$-isomorphism. 
There are some technical difficulties in proceeding this story
precisely. See \cite[subsection 4E.]{hk:fukayadeform}.

 \subsection{The category $\Mo(P)$ of weighted Morse homotopy}\label{ssec:Mo}
\label{ssec:Mo}

If we start with a toric manifold $X$ and set $\check{M}$ as
the complement of the toric divisors, 
we obtain $M$ as a torus fibration over the interior $B$ of the
dual polytope $P$. 
As we discuss in the next section, 
from the homological mirror symmetry viewpoint, 
what we should discuss is not $Fuk(M)$ but
a kind of Fukaya category $Fuk(\bar{M})$ of a torus fibration
over $P=\bar{B}$.
As an intermediate step, we consider
the category $\Mo(P)$ of weighted Morse homotopy for
the dual polytope $P$.
This $\Mo(P)$ is a generalization of the weighted Fukaya-Oh category
given in \cite{KoSo:torus} to 
the case where the base manifold has boundaries and 
critical points may be degenerate. 

Although we need only Lagrangians of constant slopes 
(which we call affine Lagrangians) in this paper, 
we formulate $\Mo(P)$ as general as possible in the following description, 
as the framework works in more general toric cases. 
In a subsequent paper \cite{FK:F1}, 
we use this formulation again to extend the result to the case of the 
Hirzebruch surface $\mathbb{F}_1$, 
where Lagrangians defined by rational functions appear.

\vspace*{0.3cm}

\noindent
{\bf Definition of $\Mo(P)$}\quad 
The definition is as follows. 
The objects of $\Mo(P)$ are 
Lagrangian sections of $\pi:M\to B$ 
satisfying certain conditions (see for instance \cite[Definition 3.1]{chan09}). 
We extend each Lagrangian section on $B$ to that on $\bar{B}$ smoothly. 
We say that two objects $L, L'$ intersect {\em cleanly}
if there exists an open set $\widetilde{B} \subset \R^n$
such that $\bar{B}\subset\widetilde{B}$
and $L, L'$ over $B$ can be extended 
to graphs of smooth sections 
over $\widetilde{B}$ so that they intersect cleanly. 
We assume that any two objects $L,L'$ intersect cleanly. 

For each $L$, we take a function $f_L$ on $\widetilde{B}$
so that $L$ is the graph of $df_L$. 
For a given ordered pair $(L,L')$,
we assign a grading $|V|$ for each connected component $V$
of the intersection $\pi(L\cap L')$ in $P=\bar{B}$
as the dimension of the stable manifold $S_v\subset \widetilde{B}$
of the gradient vector field $-\grad(f_L-f_{L'})$ with a point $v\in V$. 
This does not depend on the choice of the point $v\in V$. 
The space $\Mo(P)(L,L')$ of morphisms 
is then set to be the $\Z$-graded vector space
spanned by the connected components $V$ of $\pi(L\cap L')\in P$ 
{\em such that there exists a point $v\in V$ which is an interior point 
of $S_v\cap P\subset S_v$}. 
\footnote{We consider the Morse cohomology degree instead of the Morse homology degree.}

\def \GT{\mathcal{GT}}
\def \HGT{\mathcal{HGT}}

Now let us consider an $(l+1)$-tuple $(L_1,\dots L_{l+1})$, $l\ge 2$,
and take a generator $V_{i(i+1)}\in\Mo(P)(L_i,L_{i+1})$ 
for each $i$ and $V_{1(l+1)} \in \Mo(P)(L_1,L_{l+1})$. 
We denote by 
$$\GT(v_{12},\dots,v_{l(l+1)};v_{1(l+1)})$$
the set of gradient trees 
starting at $v_{12},\dots, v_{l(l+1)}$,
where $v_{i(i+1)}\in V_{i(i+1)}$, 
and ending at $v_{1(l+1)}\in V_{1(l+1)}$.
Here, a gradient tree 
$\gamma\in\GT(v_{12},\dots,v_{l(l+1)};v_{1(l+1)})$ is a
continuous map $\gamma:T\to P$ with a rooted trivalent $l$-tree $T$.
Regarding $T$ as a planar tree, the leaf external vertices 
from the left to the right are mapped to 
$v_{12},\dots, v_{l(l+1)}$, and
the root external vertex is mapped to $v_{1(l+1)}$ by $\gamma$. 
Furthermore, for each edge $e$ of $T$, the restriction
$\gamma|_{e}$ is a gradient trajectory of the corresponding gradient 
vector field. See \cite{KoSo:torus}. 
We then denote
\[
\GT(V_{12},\dots,V_{l(l+1)};V_{1(l+1)})
:=\bigcup_{(v_{12},\dots,v_{l(l+1)};v_{1(l+1)})
  \in V_{12}\times\cdots\times V_{l(l+1)}\times V_{1(l+1)}}
 \GT(v_{12},\dots,v_{l(l+1)};v_{1(l+1)}) . 
\]
We say that
{\em two gradient trees
$\gamma,\gamma'\in\GT(V_{12},\dots,V_{l(l+1)};V_{1(l+1)})$
are $C^\infty$-homotopic to each other} 
if $\gamma:T\to P$ is homotopic to $\gamma':T\to P$ so that 
$\gamma|_e$ is $C^\infty$-homotopic to $\gamma'|_e$ for
each edge of $T$. 
We further denote 
\[
\HGT(V_{12},\dots,V_{l(l+1)};V_{1(l+1)})
:=\{[\gamma]\ |\ \gamma\in\GT(V_{12},\dots,V_{l(l+1)};V_{1(l+1)})\}, 
\]
where $[\gamma]$ is the $C^\infty$-homotopy class of $\gamma$.

We in particular consider the case 
where $|V_{1(l+1)}|=|V_{12}|+\cdots + |V_{l(l+1)}|+2-l$,
then assume that $\HGT(V_{12},\dots,V_{l(l+1)};V_{v_{1(l+1)}})$
is a finite set.
For this we need that, in the case of the trivalent tree with one interior vertex 
for example, i) the functions assigned to each edge, which are of the form $f_L - f_{L'}$, 
are (Bott-)Morse-Smale, and ii) the (un)stable manifolds of the critical points of 
the above mentioned functions intersect transversely.
It is well-known that such transversality can always be achieved by 
a small perturbation of the functions. 
See \cite{FO} for the precise formulation.
As we shall see in the next section, 
this transversality condition is satisfied for our specific choices of Lagrangians 
and we do not further discuss this here.

For each element in $\gamma\in\GT(V_{12},\dots,V_{l(l+1)};V_{v_{1(l+1)}})$, 
we can assign the weight $e^{-A(\gamma)}$ 
where $A(\gamma)\in [0,\infty]$ is the symplectic area of 
the piecewise smooth disk in $\pi^{-1}(\gamma(T))$ 
as is done in Kontsevich-Soibelman \cite{KoSo:torus}. 
This weight is invariant with respect to a $C^\infty$-homotopy. 
Then, we define a multilinear product 
\[
m_l: \Mo(P)(L_1,L_2)\otimes \Mo(P)(L_2,L_3)\otimes
\cdots\otimes \Mo(P)(L_l,L_{l+1})\to \Mo(P)(L_1,L_{l+1})
\]
of degree $2-l$ by
\begin{equation}\label{m_l}
 m_l(V_{12}, \dots , V_{l(l+1)})=
\sum_{V_{1(l+1)}}\, \ 
\sum_{[\gamma]\in\HGT(V_{12},\dots,V_{l(l+1)};V_{1(l+1)})}
\pm e^{-A(\gamma)} V_{1(l+1)}, 
\end{equation}
where $V_{1(l+1)}$ are the bases of
$\Mo(P)(L_1,L_{n+1})$ of degree $|V_{12}|+\cdots + |V_{l(l+1)}|+2-l$.
We expect that $\{ m_l \}_{l \geq 2}$ forms a minimal $A_\infty$-structure 
with the sign $\pm$ given by the formula of 
the homological perturbation lemma.
See for example \cite[subsection 2.2]{hk:hpt-hms}.

We do not prove the $A_\infty$-relations of $m_l$'s in general for the following reason.
For a given set $\cE$ of objects of $\Mo(P)$, 
we denote by $\Mo_\cE(P)$ the full subcategory of $\Mo(P)$
consisting of objects in $\cE$, 
and consider only $\Mo_\cE(P)$ instead of the whole $\Mo(P)$.
In our examples, we can find a strongly exceptional collection $\cE$ 
in $\Tri (\Mo_\cE(P))$, 
which implies that 
all higher $m_l$'s with $l \geq 3$ vanish in $\Mo_\cE(P)$ for the degree reason.
As we shall explain later, 
we take $\cE$ to be the set of 
Lagrangian sections $L_a$ corresponding to $\cO(a)$ for the $X=\C P^n$ case, 
and take $\cE$ to be their products for the $X=\C P^m \times \C P^n$ case.

Establishing the whole $\Mo(P)$ for general $P$ requires 
much more work (see \cite{FO} 
for the case of closed manifolds).
We shall carry this out in a forthcoming paper.

\vspace*{0.3cm}

\noindent
{\bf A strong minimality assumption}\quad
For each pair $(L,L')$,
the differentials of a Morse-Bott version of the Floer complex
on $\Mo(P)(L,L')$ should be trivial in order that $\{ m_l \}_{l \geq 2}$ forms 
a minimal $A_\infty$-structure.
In our construction, we rather impose a stronger assumption
for the class of objects as follows. 
{\em For any pair $(L,L')$ and any two distinct elements 
of the basis $V,W\in\Mo(P)(L,L')$, 
there does not exist any gradient flow starting at a point in $V$ and
ending at a point in $W$. }

\vspace*{0.3cm}

\noindent
{\bf The identity morphism}\quad
For each $L\in\Mo(P)$, the space $\Mo(P)(L,L)$ of morphisms 
is generated by $P$ itself which is of degree zero. 
When the above $\Mo(P)$ is well-defined and forms
a minimal $A_\infty$-category, 
we believe that $P$ is the strict unit. 
We see that $P\in\Mo(P)(L,L)$ forms at least the identity morphism
with respect to $m_2$ under the strong minimality assumption above. 
In order to show that it is the strict unit, 
we need to show that $\Mo(P)$ is obtained
by applying homological perturbation theory to a DG category. 
Thus, that $\Mo(P)$ forms a minimal $A_\infty$-category
and that $\Mo(P)$ is strictly unital should be shown
at the same time.

\vspace*{0.3cm}

\noindent
{\bf More explicit expression}\quad 
For each $L$, let us choose a local expression
$s:B\to TB$ as we did for $\cF$ or $\cV$.
(A different choice leads to an isomorphic object. )
This enables us to assign each generator $V$ of a morphism space
a $\Z^n$-grading (which is different from the grading $|V|$ above). 
For instance, for lifts $s_a:B\to TB$ and $s_b:B\to TB$ of $L_a$ and $L_b$, 
consider $s_{b,I}:B\to TB$ defined by 
\[
 y^j=s^j_{b,I}=s^j_b -2\pi i_j, 
\]
where $I=(i_1,\dots, i_n)\in\Z^n$. 
Denote by $Mo_I(P)(s_a,s_b)$ the space generated by 
the generators of $\Mo(P)(s_a,s_b)$
which are included in the image of the intersection 
$\graph(s_a)\cap \graph(s_{b;I})$ by $TB\to B$. 
Then, we have the decomposition
\[
 \Mo(P)(s_a,s_b)=\coprod_{I\in \Z^n}\Mo_I(P)(s_a,s_b) . 
\]
In this way, each generator of $\Mo(P)(s_a,s_b)$ is assigned 
a $\Z^n$-grading. 
The multilinear product (\ref{m_l}) preserves these gradings 
in order for the corresponding gradient trees to be well-defined. 
(See \cite[subsection 4E.]{hk:fukayadeform}).

\vspace*{0.3cm}

\noindent
{\bf Connection to $\DG(X)$}\quad 
We end with this subsection by explaining
why we expect this $\Mo(P)$ to be a candidate
of the category on the mirror dual of $X$. 
We start with the DG category $\DG(X)$ of
holomorphic line bundles on $X$, 
as is constructed explicitly in subsection \ref{ssec:DGcpn}, 
and remove the toric divisors of $X$ to obtain $\check{M}$. 
Then, $\DG(X)$ should be regarded as a subcategory, 
which we denote by $\cV'(M)$, 
of $\cV(\check{M})$. 
In particular, 
$\cV'(M)$ is not full
since the smoothness condition at the removed toric divisors
is imposed in $\cV(\check{M})$. 
The cohomologies of the morphism spaces in $\cV'(M)\simeq\DG(X)$
then differs from those in $\cV(\check{M})$. 
They can be larger than those in $\cV(\check{M})$
though the morphism spaces in $\DG(X)$ are smaller at the cochain level.

In the original set-up where $M$ and $\check{M}$ are supposed to be
compact Calabi-Yau manifolds, the cohomologies of the morphism spaces
in $\cV(\check{M})$ are in one-to-one correspondence with
those in $\cF(M)$ since $\cF(M)$ and $\cV(\check{M})$ are
isomorphic DG-categories (subsection \ref{ssec:equivalent}).
Furthermore, the cohomologies of the morphism spaces $\cF(M)(s_a,s_b)$ 
are isomorphic to $Fuk(M)(s_a,s_b)$
at least when $s_a$ and $s_b$ define Lagrangian sections
$L_a$ and $L_b$ which are transversal to each other. 
Namely, the cohomologies $H(\cV(s_a,s_b))$ are spanned by
bases which are associated with 
connected components of $L_a\cap L_b$. 
We would like to keep this relation even when
$\check{M}$ is noncompact.
Then, if the smoothness condition
at the removed toric divisors 
produces additional generators in 
$H(\DG(X)(s_a,s_b))$ from $H(\cV(\check{M})(s_a,s_b))$, 
we would like to enlarge $M$ so that
there exist the corresponding additional connected
components of the intersections of the Lagrangians. 
Our feeling is that it seems to go well
if we add the boundary of $M$ and consider a kind of
Fukaya category $Fuk(\bar{M})$ or
the corresponding category $\Mo(P=\bar{B})$
of weighted Morse homotopy. 
A candidate is the category $\Mo(P)$ we defined in this subsection. 
Actually, we can consider a DG category $\cF'(M)$ which is canonically isomorphic to $\cV'(\check{M})$ 
just in a similar way as in subsection \ref{ssec:DG_M}. 
(We may just replace $d\bar{z}_i$ by $dx_i$. ) 
Then, 
we expect to have a sequence of $A_\infty$-equivalences 
\[
 \Mo(P)\simeq\cF'(M)\simeq\cV'(\check{M})\simeq\DG(X) 
\]
or the corresponding derived equivalence $\Tri(\Mo(P))\simeq\Tri(DG(X))\simeq D^b(coh(X))$. 
However, as mentioned in subsection \ref{ssec:compare}, 
there are some technical difficulties, even in the original setting, to show the $A_\infty$-equivalence $\Mo(B)\to\cF$.  
Theremore, before discussing the above equivalences in a general set-up, 
in this paper 
we explicitly proceed with this story successfully
for $X$ the projective spaces and their products 
in the next section. 
More precisely, we consider full subcategories $\Mo_\cE(P)\subset\Mo(P)$ and $\cV'_\cE(\check{M})\subset\cV'(\check{M})$ 
and an $A_\infty$-equivalence $\Mo_\cE(P)\simeq\cV'_\cE(\check{M})$, which 
induces the derived equivalence $\Tri(\Mo(P))\simeq D^b(coh(X))$. 
Note that we show the above $A_\infty$-equivalence directly and skip the intermidiate category $\cF'_\cE(M)$ there. 
Once we are convinced that $\Mo(P)$ is the correct notion, we would like to 
construct the corresponding Fukaya category $Fuk(\bar{M})$ in the future. 

 \section{Homological mirror symmetry of $\C P^n$}

 In this section, we discuss a version of homological mirror symmetry 
of $\C P^n$ as the complex side 
by explicitly proceeding with the story described in the last subsection.
In subsection \ref{ssec:DGcpn}, we construct the DG category
$\DG(\C P^n)$ of holomorphic line bundles on $\C P^n$, 
and recall the structure of its cohomologies. 
Then, we discuss the homological mirror symmetry for $\C P^n$
in subsection \ref{ssec:hms-cpn}. 
We extend the story to $\C P^m\times\C P^n$ in subsection \ref{ssec:hms-prod}.

\subsection{DG category $\DG(\C P^n)$ of line bundles over $\C P^n$}
\label{ssec:DGcpn}

We first construct the DG category $\DG(\C P^n)$ consisting of
holomorphic line bundles $\cO(a)$, $a\in\Z$. 
The space $\DG(\C P^n)(\cO(a),\cO(b))$ of morphisms is 
defined as the Dolbeault resolution of $\Gamma(\cO(a),\cO(b))$. 
Namely, it is the graded vector space,
each graded piece of which is given by 
\[
 \DG^r(\C P^n)(\cO(a),\cO(b)):=\Gamma(\cO(a),\cO(b))\otimes\Omega^{0,r}(\C P^n) 
\]
with $\Gamma(\cO(a),\cO(b))$ being the space of smooth 
bundle morphisms from $\cO(a)$ to $\cO(b)$. 
The composition of morphisms is defined in a similar way as
that in $\cV(\check{M})$ in subsection \ref{ssec:DG_cM}. 
Each $\cO(a)$ is associated with the connection $D_a$,
which is expressed locally as 
\[
 D_a = d -a
\frac{\bar{w}_1dw_1+\cdots + \bar{w}_ndw_n}{1+\bar{w}_1w_1+\cdots +\bar{w}_nw_n}
\]
on $U=U_0$ (eq.(\ref{conn-cpn})),
and the differential 
\[
 d_{ab}: \DG^r(\C P^n)(\cO(a),\cO(b))\to \DG^{r+1}(\C P^n)(\cO(a),\cO(b))
\]
is defined by
\[
 d_{ab}(\ti\psi)
 := 2\left(D_b^{0,1}\ti\psi- (-1)^r\ti\psi D_a^{0,1}\right) . 
\]
This differential satisfies the Leibniz rule with respect to
the composition. Thus, $\DG(\C P^n)$ is a DG category. 

The generators of $H^0(\DG(\C P^n))(\cO(a),\cO(a+1))$ are given by
\begin{equation}\label{gen-cpn}
 1,\ w_1 ,\ w_2,\ \dots,\ w_n 
\end{equation}
locally on $U$. 
These generate $H^0(\DG(\C P^n))(\cO(a),\cO(b))$ as products of functions, 
so $H^0(\DG(\C P^n))(\cO(a),\cO(b))$ is represented by 
polynomials in $(w_1,\dots, w_n)$ of degree equal to or less than $b-a$. 
In particular, $H^0(\DG(\C P^n))(\cO(a),\cO(b))=0$ for $a>b$.  
It is known by \cite{beilinson78} that 
$\cE:=(\cO(q),\dots,\cO(q+n))$
forms a full strongly exceptional collection of $D^b(coh(\C P^n))$
for each $q\in\Z$. 
That $\cE$ forms a strongly exceptional collection 
means
\begin{equation*}
 \begin{split}
   & H^0(DG(\C P^n))(\cO(a),\cO(a))\simeq\C , \\
   & H^0(DG(\C P^n))(\cO(a),\cO(b))=0,\qquad a>b \\
   & H^r(DG(\C P^n))(\cO(a),\cO(b))=0,\qquad r\ne 0
  \end{split}
\end{equation*}
for any $a,b=\{q,q+1,\dots, q+n\}$. 
Let $\DG_\cE(\C P^n)$ be the full DG subcategory of $DG(\C P^n)$
consisting of $\cE$.
Then the strongly exceptional collection $\cE$ is full means that 
it generates $D^b(coh(\C P^n))$ in the sense that 
\[
 \Tri(\DG_{\cE}(\C P^n))\simeq D^b(coh(\C P^n)),  
\]
where $\Tri$ is the Bondal-Kapranov construction \cite{BK:enhanced}. 

Note that $H^{n}(\DG(\C P^n))(\cO(q+n+1),\cO(q))\ne 0$;  
it includes an element represented by 
\[
 \frac{d\bar{w}_1\cdots d\bar{w}_n}{(1+w_1\bar{w}_1+\cdots + w_n\bar{w}_n)^2} . 
\]

 \subsection{Homological mirror symmetry of $\C P^n$}
\label{ssec:hms-cpn}

First, we identify the DG category $\DG(\C P^n)$
with a (non-full) subcategory $\cV'$ of
the DG category $\cV=\cV(\check{M})$ 
consisting of the same objects $\cO(a)$, $a\in\Z$, 
where
\[
\check{M}=\C P^n \backslash
 \{[t_0:t_1:\cdots : t_n]\ |\ t_0\cdot t_1\cdots t_n=0\} . 
\]
For a given morphism
$\ti\psi\in\DG^0(\C P^n)(\cO(a),\cO(b))$, 
we express it locally on $\cU$ (see subsection \ref{ssec:tt}), and remove
the origin (corresponding to $t_0=0$). 
We send this to $\cV'(\cO(a),\cO(b))$ using (\ref{Psi_a}): 
\[
 \ti\psi\mapsto \psi:=\Psi_b^{-1}\circ\ti\psi\circ\Psi_a . 
\]
Clearly, this map is compatible with
the differentials and the compositions in both sides. 
In this way, we obtain a functor
\[
 \cI:\DG(\C P^n)\to \cV
\]
of DG-categories. 
We see that $\cI$ is faithful.
However, $\cI$ is not full since 
$\ti\psi$ is smooth at the points
$\{[t_0:t_1:\cdots : t_n]\ |\ t_0\cdot t_1\cdots t_n=0\}$. 
Thus, the image
\[
 \cV':=\cI(\DG(\C P^n))
\]
is a non-full DG subcategory of $\cV$.

The local expression for morphisms are transformed by $\cI$ as follows. 
By $w_i=e^{x_i+\ii y_i}=r_i e^{\ii y_i}$, we Fourier-expand $\ti\psi$ as
\[
 \ti\psi=\sum_{I\in\Z^n}\ti\psi_I(r) e^{\ii I\check{y}} ,
 \qquad r:=(r_1,\dots, r_n). 
\]
Now, recall (\ref{dual-x}) and then 
\[
x^1+\cdots + x^n
= \frac{2(e^{2x_1}+\cdots + e^{2x_n})}{1+e^{2x_1}+\cdots + e^{2x_n}}
= 2-\frac{2}{1+e^{2x_1}+\cdots + e^{2x_n}} , 
\]
so we have
\[
\Psi_a=\left( 1+e^{2x_1}+\cdots + e^{2x_n}\right)^{\frac{a}{2}}
 = \left(\frac{2}{2-x^1-x^2-\cdots -x^n}\right)^{\frac{a}{2}} , 
\]
and 
\[
 r_i=e^{x_i}= \left(x^i\frac{1+e^{2x_1}+\cdots + e^{2x_n}}{2}\right)^{\ov{2}}
 = \left(\frac{x^i}{2-x^1-x^2-\cdots -x^n} \right)^{\ov{2}} . 
\]
Then, $\psi=\cI(\ti\psi)$ turns out to be 
\[
\sum_{I\in\Z^n}\ti\psi_I(r)\left(\frac{2-x^1-x^2-\cdots -x^n}{2}\right)^{\frac{b-a}{2}}e^{\ii I \check{y}}
=\sum_{I\in\Z^n}\psi(x) e^{\ii I\check{y}} , 
\]
so Fourier-componentwisely we have the transformation 
\[
 \psi_I(x)=\ti\psi_I(r(x))\left(\frac{2-x^1-x^2-\cdots -x^n}{2}\right)^{\frac{b-a}{2}} . 
\]
We can bring the generators (\ref{gen-cpn})
of $H^0(\DG(\C P^n)(\cO(a),\cO(a+1)))$ over $\C$ 
to those of $H^0(\cV')(\cO(a),\cO(a+1))$, which are given by 
\begin{equation}\label{base1}
 \left[\sqrt{\frac{2-x^1-x^2-\cdots - x^n}{2}}\right]
\end{equation}
and 
\begin{equation}\label{base2}
 \left[\sqrt{\frac{x^1}{2}} e^{\ii y_1}\right] ,\
 \left[\sqrt{\frac{x^2}{2}} e^{\ii y_2}\right] ,\
 \dots, \left[\sqrt{\frac{x^n}{2}} e^{\ii y_n}\right]  . 
\end{equation}
The above bases (\ref{base1}) and (\ref{base2}) 
generate the whole space $H^0(\cV')(\cO(a),\cO(b))$ 
as products of these functions. 
Explicitly, the bases ${\bf e}_{ab;I}$, $I=(i_1,\dots, i_n)$,
of the vector space 
$H^0(\cV')(\cO(a),\cO(b))$ are 
\begin{equation}\label{e_abI}
 {\bf e}_{ab;I}= 
c_{ab;I}\cdot\left(\sqrt{\frac{2-x^1-x^2-\cdots - x^n}{2}}\right)^{b-a-|I|}
 \left(\sqrt{\frac{x^1}{2}}e^{\ii y^1}\right)^{i_1}
 \cdots
  \left(\sqrt{\frac{x^n}{2}}e^{\ii y^n}\right)^{i_n} , 
\end{equation}
where $i_1\ge 0,\dots, i_n\ge 0$ and $|I|:=i_1+\cdots + i_n\le b-a$,
and we attach $c_{ab;I}$ 
so that $\max_{x\in P} |{\bf e}_{ab;I}(x)|=1$. 
Note that this is valid for $a=b$, too, where
we only have $I=(0,\dots,0)=:0$ and $e_{aa;0}$ is the
identity element; $e_{aa;0}(x)=1$ for any $x\in P$. 

Since all the exponents in (\ref{e_abI}) are non-negative, 
we see that each $e_{ab;I}$ extends to a continuous function on $P$. 
By direct calculations, we have the following lemma. 
\begin{lem}\label{lem:2}
For a fixed $a<b$ and ${\bf e}_{ab;I}\in H^0(\cV')(\cO(a),\cO(b))$, 
the set 
\[
 \{x\in P\ | |{\bf e}_{ab;I}(x)|=1\}
\]
consists of a point
\[
 v_{ab;I}:=\left(\frac{2 i_1}{b-a},\dots, \frac{2 i_n}{b-a} \right), 
\]
which is the intersection $V_{ab;I}\subset \pi(L_a\cap L_b)$
with label $I$. 
This correspondence then gives a quasi-isomorphism
\[
\iota:\Mo(P)(L_a,L_b)\to \cV'(\cO(a),\cO(b))
\]
of cochain complexes. 
\qed\end{lem}
For each $a<b$ and $I$, we later employ a function $f_{ab;I}$ on $P$ defined uniquely by 
\begin{equation}\label{fab}
 \sum_{j=1}^n (s_a^j-s_b^j+2\pi i_j)dx_j = df_{ab;I} , \qquad f_{ab;I}(v_{ab;I})=0 .  
\end{equation}
\begin{rem}\label{fab}
In the original set-up where $B$ is compact, 
this $f_{ab;I}$ is a Morse function, where $v_{ab;I}$ is the critical point
of degree zero. 
However, now in our case, $v_{ab;I}$ may be at the boundary $\partial(P)$. 
Even if we extend $P$ to $\ti{B}$ naturally, $v_{ab;I}\in\partial(P)$
may not be a critical point since the symplectic form on $M$ diverges at
the boundary. 
\end{rem}

For each $a$, the space $\Mo(P)(L_a,L_a)$ is generated by $P$. 
The two conditions 
\[
 \max_{x\in P} |e_{aa;0}(x)| =1,\qquad \{x\in P\ |\ |e_{aa;0}(x)|=1 \}= P
\]
are clearly satisfied. We define a quasi-isomorphism 
$\iota:\Mo(P)(L_a,L_a)\to\cV'(\cO(a),\cO(a))$ by
$\iota(P)=e_{aa;0}$. 

For $a>b$, both the space $Mo(P)(L_a,L_b)$ and the cohomology
$H(\cV'(\cO(a),\cO(b)))$ are trivial.
Thus, the zero map $\iota:\Mo(P)(L_a,L_b)\to\cV'(\cO(a),\cO(b))$
is a quasi-isomorphism.

Now, let us fix $q\in\Z$ and consider
$\cE:=(\cO(q),\cO(q+1),\dots,\cO(q+n))$. We denote the
corresponding full subcategories by 
$\DG_{\cE}(\C P^n)\subset\DG(\C P^n)$, 
$\cV'_{\cE}\subset\cV'$ and $\Mo_{\cE}(P)\subset \Mo(P)$. 
It is known that $\cE$ (with any $q$) forms
a full strongly exceptional collection
in $\Tri(\DG_{\cE}(\C P^n))\simeq D^b(coh(\C P^n))$ 
\cite{bondal:asso}. 
Recall that an $A_\infty$-equivalence is an $A_\infty$-functor
which induces a category equivalence
on the corresponding cohomology categories. 
\begin{thm}\label{thm:cpn}
For each $q\in\Z$, the quasi-isomorphisms
\[
 \iota:\Mo(P)(L_a,L_b)\to\cV'(\cO(a),\cO(b))
\]
with $a,b\in \{q,...,q+n\}$ 
extend to a linear $A_\infty$-equivalence
\[
 \iota: \Mo_\cE(P)\overset{\sim}{\to}\cV'_\cE .
\]
\end{thm}
Now, 
we have the DG isomorphism $\DG_\cE(\C P^n)\simeq\cI(DG_\cE(\C P^n))=\cV'_\cE$. 
Since a DG functor is a linear $A_\infty$-functor, we immediately obtain the
following.
\begin{cor}\label{cor:cpn1}
One has a linear $A_\infty$-equivalence 
\[
 \Mo_{\cE}(P)\to \DG_\cE(\C P^n) . 
\]
\qed\end{cor}
\begin{cor}\label{cor:cpn2}
One has an equivalence of triangulated categories
\[
 \Tri(\Mo_{\cE}(P))\simeq D^b(coh(\C P^n)) . 
\]
\qed\end{cor}

In the rest of this subsection, we show Theorem \ref{thm:cpn} 
by computing the structure of $\Mo_{\cE}(P)$. 
We already see that nontrivial morphisms in $\Mo(P)$ are of degree zero only. 
This implies, by degree counting, that the higher $A_\infty$-products of
$\Mo(P)$ are trivial. 
Thus, what remains to show 
the theorem is to construct
the product $m_2$ and show
the compatibility of the products with respect to $\iota$.  
\begin{lem}\label{lem:3}
For $a<b<c$ and bases 
$V_{ab;I_{ab}}\in\Mo(P)(L_a,L_b)$, $V_{bc;I_{bc}}\in\Mo(P)(L_b,L_c)$, 
we have 
\begin{equation}\label{compati-lem3}
 \iota\, m_2(V_{ab;I_{ab}},V_{bc;I_{bc}})
 = {\bf e}_{ab;I_{ab}}\cdot {\bf e}_{bc;I_{bc}} . 
\end{equation}
\end{lem}
\begin{pf}
Recall that each base consists of a point; $V_{ab;I_{ab}}=\{v_{ab;I_{ab}}\}$ 
and so on. 
We take the function $f_{ab;I_{ab}}$ defined by (\ref{fab}). 
Since its gradient vector field is of the form 
\[
 -\grad(f_{ab;I_{ab}})= 2\pi \frac{(b-a)}{2}
 \left((x^1-i_{ab;1})\fpartial{x^1}+\cdots
 + (x^n-i_{ab;n})\fpartial{x^n}\right) , 
\]
its gradient trajectories starting from $v_{ab;I_{ab}}$ go straight. 
Similarly, gradient trajectories of $-\grad(f_b-f_c)$
starting from $v_{bc;I_{bc}}$ go straight. 
On the other hand, the only gradient trajectory of $-\grad(f_a-f_c)$ ending at
$v_{ac;I_{ac}}$, $I_{ac}:=I_{ab}+I_{bc}$
is the one staying at $v_{ac;I_{ac}}$ 
since it is of degree zero.
This means that these three gradient trajectories should meet at
$v_{ac;I_{ac}}$. 
Thus we obtained the gradient tree $\gamma$ defining the product
$m_2(V_{ab;I_{ab}},V_{bc;I_{bc}})$ explicitly. 
(The result is that $v_{ac;I_{ac}}$ sits on the straight line segment $v_{ab;I_{ab}}v_{bc;I_{bc}}$ 
in all cases. )
Now, $A(\gamma)$ turns out to be 
\[
 A(\gamma)= \frac{1}{2\pi} f_{ab;I_{ab}}(v_{ac;I_{ac}})+ \frac{1}{2\pi} f_{bc;I_{bc}}(v_{ac;I_{ac}}) . 
\]
Here, $f_{ab;I_{ab}}(v_{ac;I_{ac}})/2\pi$ is the symplectic area
of the triangle disk whose edges belong to
$s_a(\gamma(T))$, $s_b(\gamma(T))$ and $\pi^{-1}(v_{ac;I_{ac}})$. 
Similarly, $f_{bc;I_{bc}}(v_{ac;I_{ac}})/2\pi$ is 
the symplectic area of the corresponding triangle disk. 
We thus obtain the weight $+e^{-A(\gamma)}$. 

Next, we look at the product in $\cV'$ side. 
We can express the bases as 
${\bf e}_{ab;I_{ab}}=e^{(1/2\pi) f_{ab;I_{ab}}}\cdot e^{\ii I_{ab} \check{y}}$ 
and 
${\bf e}_{bc;I_{bc}}=e^{(1/2\pi) f_{bc;I_{bc}}}\cdot e^{\ii I_{bc} \check{y}}$. 
We have
\[
  {\bf e}_{ab;I_{ab}}\cdot {\bf e}_{bc;I_{bc}}
  = e^{\frac{1}{2\pi}\left(f_{ab;I_{ab}}+f_{bc;I_{bc}}\right)}\cdot e^{\ii I_{ac} \check{y}} . 
\]
Since this is the product of the zero-th cohomologies, the result is also 
a closed morphism in $\cV'$. 
Hence, the right hand side is proportional to ${\bf e}_{ac;I_{ac}}$, 
whose absolute value takes the maximal value at $v_{ac;I_{ac}}$. 
Namely, we have
\[
  {\bf e}_{ab;I_{ab}}\cdot {\bf e}_{bc;I_{bc}}
  = e^{\frac{1}{2\pi}\left(f_{ab;I_{ab}}(v_{ac;I_{ac}})+f_{bc;I_{bc}}(v_{ac;I_{ac}})\right)} \cdot {\bf e}_{ac;I_{ac}} . 
\]
This shows that the compatibility (\ref{compati-lem3}) holds true. 
\qed\end{pf}

We need to show the compatibility (\ref{compati-lem3}) for 
any $a\le b\le c$. 
If $a=b$, then $V_{ab;I_{ab}}=P$. If $b=c$, then $V_{bc;I_{bc}}=P$. 
Now, we see that $\Mo_\cE(P)$ satisfies
the strong minimality assumption in subsection \ref{ssec:Mo},
which implies that $P$ forms the identity morphism in $\Mo_\cE(P)$. 
Since we already know that $\iota(P)$ is the identify morphism
in $\cV'_\cE$, the compatibility (\ref{compati-lem3}) follows and 
the proof of 
Theorem \ref{thm:cpn} is completed. \qed
\begin{rem}\label{rem:hpt}
The product $m_2$ and the linear $A_\infty$-equivalence $\iota$ 
can be induced by applying homological perturbation theory
to $DG(\C P^n)$ in a suitable way. 
As the higher $A_\infty$-products of $\Mo(P)$ are trivial,
the induced $A_\infty$-equivalence turns out to be linear 
by degree counting 
since nontrivial cohomologies of morphisms are of degree zero only. 
\end{rem}

As a biproduct of the proof, we see that $\Mo_{\cE}(P)$ has
the following properties.
\begin{prop}\label{prop:boundary}
For any $L_a, L_b\in\Mo_{\cE}(P)$ such that $L_a\ne L_b$,
$V_{ab}=\pi(L_a\cap L_b)$ belongs to the boundary $\partial(P)$. 

For given bases
$V_{ab}\in\Mo_{\cE}(L_a,L_b)$ and
$V_{bc}\in\Mo_{\cE}(L_b,L_c)$, 
the image $\gamma(T)$ by any gradient tree
$\gamma\in\GT(V_{ab},V_{bc};V_{ac})$ 
belongs to the boundary $\partial(P)$
unless $L_a=L_b=L_c$. 
\qed\end{prop}
\begin{rem}
If $L_a=L_b=L_c$, then $V_{ab}=P$, $V_{bc}=P$ and $V_{ac}=P$. 
Then $\gamma\in\GT(V_{ab},V_{bc};V_{ac})$ is a constant map
to a point in $P$.
If $L_a=L_b\ne L_c$, then $V_{ab}=P$ and $V_{bc}=\{v_{bc}\}=V_{ac}$.
Then $\gamma\in\GT(V_{ab},V_{bc};V_{ac})$ is the constant map 
to the point $v_{bc}\in\partial(P)$. 
Similarly, if $L_a\ne L_b=L_c$, 
then $\gamma\in\GT(V_{ab},V_{bc};V_{ac})$ is the constant map
to the point $v_{ab}\in\partial(P)$. 
\end{rem}
We expect that for many other toric Fano manifolds $X$ and
(strongly) exceptional collections $\cE$, 
$\Mo_\cE(P)$ may satisfy these properties.

We also believe that there exists an $A_\infty$-equivalence
\[
 \Mo(P)\to \DG(\C P^n)
\]
between the whole categories. 
However, it is not easy to show directly
that the whole category $Mo(P)$ is well-defined as an $A_\infty$-category
since there are infinitely many gradient trees for which
we should check whether our assumption holds or not. 
In particular, if $\Mo(P)$ is well-defined, it should have nontrivial higher $A_\infty$-products.

  \subsection{Homological mirror symmetry of $\C P^m\times\C P^n$}
\label{ssec:hms-prod}

In this subsection we shall see how the framework presented in 
the last subsections works for the case of the product of projective spaces. 
The point here is that we need not only transversal but clean intersections 
of Lagrangians in the symplectic side. 
That's why we included clean intersections in the definition of 
$\Mo(P)$ in subsection 4.5. 
We still do not need higher products $m_3$, $m_4$,... because 
we can pick up full strongly exceptional collections on both sides 
(see remark at the end of subsection 4.5). 

Let $X=\C P^m\times \C P^n$, and $\check{M}$ be
the complement of the toric divisors. 
For
\[
\xymatrix{
& X\ar[ld]_{p_1}\ar[rd]^{p_2} & \\
\C P^m &  & \C P^n
}
\]
we denote $\cO(a,b):=p_1^*\cO(a)\otimes p_2^*\cO(b)$.
Then $\cE:=\{ \cO(a,b) \}_{a=0,1,...,m,\, b=0,1,...,n}$ 
with the lexicographic order 
forms a strongly exceptional collection.
According to Orlov \cite{orlov:92} the semi orthogonal ordered set of 
admissible subcategories $(\mathcal{D}_0,...,\mathcal{D}_n)$, 
where $\mathcal{D}_0$ is the image of $D^b(coh(\C P^m))$ 
in $D^b(coh(\C P^m \times \C P^n))$ under the pull-back functor $p_1^\ast$ 
and $\mathcal{D}_i$'s are its twists 
along $\C P^n$, generates $D^b(coh(\C P^m \times \C P^n))$, 
which means that the collection of $\cO(a,b)$'s above is full.
We consider the DG category $\DG(X)$ of these line bundles,
and the corresponding DG-category $\cV(\check{M})$. 
Just in a similar way as in the previous subsection,
we have a DG subcategory $\cI(\DG(X))=\cV'\subset\cV=\cV(\check{M})$
so that $\DG(X)\simeq\cV'$.
Their full subcategories consisting of $\cE$ are denoted
$\DG_\cE(X)$ and $\cV'_\cE$. 

Then, the parallel statements to the case $X=\C P^n$ 
hold.
\begin{thm}\label{thm:prod}
There exists a linear $A_\infty$-equivalence
\[
 \iota: \Mo_\cE(P)\overset{\sim}{\to}\cV'_\cE
\]
such that for each generator $V\in\Mo_\cE(P)(L,L')$ 
we have 
$\max_{x\in P} |\iota(V)(x)|=1$ and 
\[
 V=\{x\in P\ |\ |\iota(V)(x)|=1\ \} . 
\]
\end{thm}
\begin{cor}\label{cor:prod1}
We have a linear $A_\infty$-equivalence
\[
 \Mo_\cE(P)\simeq \DG_\cE(\C P^m\times\C P^n) . 
\]
\qed\end{cor}
\begin{cor}\label{cor:prod2}
We have an equivalence of triangulated categories 
\[
 \Tri(\Mo_\cE(P))\simeq D^b(coh(\C P^m\times\C P^n)) . 
\]
\qed\end{cor}
\begin{prop}\label{prop:prod}
If $L\ne L'$, any generator $V\in\Mo_{\cE}(P)(L,L')$
belongs to the boundary $\partial(P)$. 

For bases
$V\in\Mo_{\cE}(L,L')$ and
$V'\in\Mo_{\cE}(L',L'')$ such that 
$L\ne L'$ and $L'\ne L''$, any gradient tree
$\gamma\in\GT(V,V';V'')$ with $V''\in\Mo_\cE(P)(L,L'')$ 
belongs to the boundary $\partial(P)$. 
\end{prop}

\noindent
{\em proof of Theorem \ref{thm:prod}}\quad 
The bases of the space $H^0(\cV')(\cO(a_1,a_2),\cO(b_1,b_2))$ are 
\[
 e_{a_1b_1;I}\otimes e_{a_2b_2;J}
\]
where $e_{a_1b_1;I}$ and $e_{a_2b_2;J}$ are
the bases of the corresponding zero-cohomology spaces of morphisms
defined in (\ref{e_abI}) 
for $\C P^m$ and $\C P^n$, respectively, so 
$I=(i_1,\dots,i_m)$ and $J=(j_1,\dots,j_n)$ run over
\[
 \begin{split}
 & i_1\ge 0,\dots, i_m\ge 0,\quad |I|\le b_1-a_1, \\
 & j_1\ge 0,\dots, j_n\ge 0,\quad |J|\le b_2-a_2 . 
 \end{split}
\]
Each base satisfies $\max_{x\in P}|(e_{a_1b_1;I}\otimes e_{a_2b_2;J})(x)|=1$. 
Let us denote by $L_{(a_1,a_2)}\in\Mo(P)$ the object corresponding to
$\cO(a_1,a_2)$. 
The base corresponding to $e_{a_1b_1;I}\otimes e_{a_2b_2;J}$ is then 
\[
 \{x\in P\ |\ |(e_{a_1b_1;I}\otimes e_{a_2b_2;J}) (x)|=1\} = 
V_{a_1b_1;I}\times V_{a_2b_2;J} \in\Mo(P)(L_{(a_1,a_2)},L_{(b_1,b_2)}) .
\]
It consists of the point $(v_{a_1b_1;I},v_{a_2b_2;J})$ 
if $a_1<b_1$ and $a_2<b_2$. 
Otherwise, we have 
\[
 V_{a_1b_1;I}\times V_{a_2b_2;J} = P_1\times \{v_{a_2b_2;J}\}
\]
for $a_1=b_1$ and $a_2<b_2$, 
\[
 V_{a_1b_1;I}\times V_{a_2b_2;J} = \{v_{a_1b_1;I}\}\times P_2
\]
for $a_1<b_1$ and $a_2=b_2$, and 
\[
 V_{a_1b_1;I}\times V_{a_2b_2;J} = P_1\times P_2 =P
\]
for $a_1=b_1$ and $a_2=b_2$,
where $P_1$ and $P_2$ are the dual polytope of
$\C P^m$ and $\C P^n$, respectively. 

For
$V_{a_1b_1;I}\times V_{a_2b_2;J}\in\Mo(P)(L_{(a_1,a_2)},L_{(b_1,b_2)})$ 
and
$V_{b_1c_1;K}\times V_{b_2c_2;L}\in\Mo(P)(L_{(b_1,b_2)},L_{(c_1,c_2)})$, 
the equation 
\[
 \iota m_2(V_{a_1b_1;I}\times V_{a_2b_2;J},V_{b_1c_1;K}\times V_{b_2c_2;L})
 = (e_{a_1b_1;I}\otimes e_{a_2b_2;J})\cdot
(e_{b_1c_1;K}\otimes e_{b_2c_2;L})  
\]
follows immediately from
looking at the structure of the gradient tree $\gamma$ 
defining 
the product $m_2(V_{a_1b_1;I}\times V_{a_2b_2;J},V_{b_1c_1;K}\times V_{b_2c_2;L})$. 
Actually, 
let us denote by $\gamma_1$ the gradient tree 
obtained as the composition of $\gamma$ with the projection $P\to P_1$.
We see that $\gamma_1$ is the gradient tree
defining the product
$m_2(V_{a_1b_1;I},V_{b_1c_1;J})$. 
Similarly, we consider $\gamma_2$.
Then, we have $A(\gamma)=A(\gamma_1)+A(\gamma_2)$, 
and we see that this is compatible with the
product 
\[
 \begin{split}
(e_{a_1b_1;I}\otimes e_{a_2b_2;J})\cdot
(e_{b_1c_1;K}\otimes e_{b_2c_2;L})
& =(e_{a_1b_1;I}\cdot e_{b_1c_1;K})\otimes (e_{a_2b_2;J}\cdot e_{b_2c_2;L}) \\
& = e^{-A(\gamma_1)}e_{a_1c_1;I+K}\otimes e^{-A(\gamma_2)}e_{a_2c_2;J+L} . 
 \end{split}
\]
This completes the proof of Theorem \ref{thm:prod}. 
\qed

\noindent
{\em Proof of Proposition \ref{prop:prod}}\quad
Each $\gamma$ is obtained from the pair $(\gamma_1,\gamma_2)$
in the proof of Theorem \ref{thm:prod}. 
In particular, the image $\gamma(T)$ of a trivalent tree $T$ by $\gamma$
is obtained as $(\gamma_1,\gamma_2)(T)\subset P_1\times P_2$.
By Proposition \ref{prop:boundary}, $\gamma_1(T)\subset \partial(P_1)$ if $V_{a_1b_1}\neq P_1$ or $V_{b_1c_1}\neq P_1$, i.e., if $a_1,b_1,c_1$ do not satisfy $a_1=b_1=c_1$. 
Similarly, 
$\gamma_2(T)\subset \partial(P_2)$ unless $a_2=b_2=c_2$. 

On the other hand, when at least
either $L_{(a_1,a_2)}\neq L_{(b_1,b_2)}$ or $L_{(b_1,b_2)}\neq L_{(c_1,c_2)}$
is satisfied, at least one of the inequalities
$a_1<b_1, b_1<c_1, a_2<b_2, b_2<c_2$ is satisfied.
Thus, at least neither $a_1=b_2=c_1$ nor $a_2=b_2=c_2$ is satisfied.
This means that at least either $\gamma_1(T)\subset \partial(P_1)$ or
$\gamma_2(T)\subset\partial(P_2)$ holds, which implies 
that $\gamma(T)\subset\partial(P)$. 
\qed

Lastly, we show more explicitly the gradient trees $\gamma$
defining the products
\[
 m_2(V_{a_1b_1;I}\times V_{a_2b_2;J},V_{b_1c_1;K}\times V_{b_2c_2;L}) . 
\]
A different point from the previous subsection is that
$V_{a_1b_1;I}\times V_{a_2b_2;J}$, $V_{b_1c_1;K}\times V_{b_2c_2;L}$
and $V_{a_1c_1;I+K}\times V_{a_2c_2;J+L}$ may not be points
even if $V_{a_1b_1;I}\times V_{a_2b_2;J}\ne P$ and
$V_{b_1c_1;K}\times V_{b_2c_2;L}\ne P$. 
There are $2^4=16$ types of the products 
depending on whether each $"\le"$ in $(a_1\le b_1\le c_1; a_2\le b_2\le c_2)$ is $"="$ or $"<"$. 
When $"="$ is included, then the corresponding morphism consists of a clean intersection 
(with nonzero dimension). 
We further divide them by the number of $"="$ as follows. 
\begin{itemize}
 \item[(0)] the case $a_1<b_1<c_1$, $a_2<b_2<c_2$. 
 \item[(1)] only one of the equations $a_1=b_1$, $b_1=c_1$, $a_2=b_2$, $b_2=c_2$ 
is satisfied. 
 \item[(2)] $a_1=b_1<c_1$ and $a_2=b_2<c_2$, or
$a_1<b_1=c_1$ and $a_2<b_2=c_2$. 
 \item[(2')] $a_1=b_1<c_1$ and $a_2<b_2=c_2$, or
$a_1<b_1=c_1$ and $a_2=b_2<c_2$. 
 \item[(2'')] $a_1=b_1=c_1$ and $a_2<b_2<c_2$, or 
$a_1<b_1<c_1$ and $a_2=b_2=c_2$. 
 \item[(3)] three of the equations $a_1=b_1$, $b_1=c_1$, $a_2=b_2$, $b_2=c_2$ 
are satisfied and the remaining one is the inequality. 
 \item[(4)] $a_1=b_1=c_1$ and $a_2=b_2=c_2$. 
\end{itemize}
The case (4) corresponds to $P \cdot P = P$ in $\Mo(P)$. 
The cases (3) and (2) correspond to
the products $P\cdot V=V$ or $V\cdot P = P$ for some $V$. 
The case (2'') corresponds to 
the product $(P_1\times V)\cdot (P_1\times W)=P_1\times (V\cdot W)$ 
or $(V\times P_2)\cdot (W\times P_2)=(V\cdot W)\times P_2$,
where $V\cdot W$ is a product in $\Mo(P_1)$ or $\Mo(P_2)$. 
Thus, the argument reduces to the one in the previous subsection. 
In the case (0),
all $V_{a_1b_1;I}\times V_{a_2b_2;J}$, 
$V_{b_1c_1;K}\times V_{b_2c_2;L}$ and
$V_{a_1c_1;I+K}\times V_{a_2c_2;J+K}$ 
consist of points. Thus, the situation is just the product of 
the ones in Lemma \ref{lem:3}. 
 
Now we discuss more carefully the cases (1) and (2'). 
In the case (1), if $a_1=b_1$, then 
\[
 V_{a_1b_1;I}\times V_{a_2b_2;J} = P_1\times \{v_{a_2b_2;J}\}  
\]
which is not a point. 
In this case, a gradient tree
$\gamma\in\GT(P_1\times \{v_{a_2b_2;J}\}, \{v_{a_1c_1;K}\}\times\{v_{b_2c_2;L}\}; \{v_{a_1c_1;K}\}\times\{v_{a_2c_2;J+L}\} )$
belongs to 
$\GT((v_{a_1c_1;K},v_{a_2b_2;J}), (v_{a_1c_1;K},v_{b_2c_2;L}); (v_{a_1c_1;K},v_{a_2c_2;J+L}))$ 
and the image $\gamma(T)$ is a straight segment 
connecting $(v_{a_1c_1;K},v_{a_2b_2;J})$ and $(v_{a_1c_1;K},v_{b_2c_2;L})$ on which
$(v_{a_1c_1;K},v_{a_2c_2;J+L})$ sits. 
Similarly, in the case (2'), if $a_1=b_1$ and $b_2=c_2$, then we have
\[
 V_{a_1b_1;I}\times V_{a_2b_2;J} = P_1\times \{v_{a_2b_2;J}\}  ,\qquad 
 V_{b_1c_1;K}\times V_{b_2c_2;L} = \{v_{b_1c_1;K}\}\times P_2 . 
\]
A gradient tree 
$\gamma\in \GT(P_1\times\{v_{a_2b_2;J}\}, \{v_{b_1c_1;K}\}\times P_2;
\{v_{b_1c_1;K}\}\times\{v_{a_2b_2;J}\})$ 
belongs to
$\GT((v_{b_1c_1;K},v_{a_2b_2;J}), (v_{b_1c_1;K},v_{a_2b_2;J});
(v_{b_1c_1;K},v_{a_2b_2;J}))$, 
and then the image $\gamma(T)$ is just the point
$(v_{b_1c_1;K},v_{a_2b_2;J})$ which is the intersection 
$(P_1\times \{v_{a_2b_2;J}\})\cap(\{v_{b_1c_1;K}\}\times P_2)$. 
Though this example is still too simple in this sense, 
this is actually an example of products of Morse homotopy 
where the gradient trees start from 
(non-transverse) clean intersections 
instead of intersection points, 
generalizing the original set-up of Morse homotopy \cite{fukayaAinfty, FO} and 
Kontsevich-Soibelman \cite{KoSo:torus}. 
More examples of clean intersections 
appear in the case of $\mathbb{F}_1$ \cite{FK:F1}.

\vspace*{0.5cm}

\noindent
{\bf ACKNOWLEDGMENTS}
\vspace*{0.2cm}

We are grateful to Akira Ishii for introducing Orlov's paper \cite{orlov:92} to us.
We also thank the referee for valuable suggestions 
which improve expositions in this paper very much.

M.~F.\ is supported by Grant-in-Aid for Scientific Research (C) (18K03269)  
of the Japan Society for the Promotion of Science. 
H.~K.\ is supported by Grant-in-Aid for Scientific Research (C) (18K03293)  
of the Japan Society for the Promotion of Science.

\end{document}